\documentclass[final,leqno]{siamltex}
\pdfoutput=1
\usepackage{amsmath}
\usepackage{graphicx}
\newtheorem{assumption}{Assumption}[section]
\newtheorem{remark}{Remark}[section]
\def\reff#1{(\ref{#1})}
\def\eps{{\varepsilon}}
\def\vphi{{\varphi}}
\def\na{{\nabla}}
\def\De{{\Delta}}
\def\de{{\delta}}
\def\Om{{\Omega}}
\def\cT{{\mathcal{T}_h}}
\def\cF{{\mathcal{F}_h}}
\def\tK{{\widetilde{K}}}
\def\G{{\Gamma}}
\def\p{{\partial}}

\def\shml{{\widetilde{H}^{-1}}}
\def\shmt{{\widetilde{H}^{-2}}}

\def\shmlt{{\widetilde{H}^{-2}(\Omega)}}
\newcommand{\norm}[2]{\left\Vert#1\right\Vert_{#2}}
\newcommand{\normL}[1]{\left\Vert#1\right\Vert_{L^2}}
\newcommand{\dual}[2]{\left\langle#1, #2\right\rangle}
\newcommand{\abs}[1]{\left\vert#1\right\vert}
\newcommand{\set}[1]{\left\{#1\right\}}

\begin{document}

\title{A posteriori error estimates for finite element approximations
of the Cahn-Hilliard equation and the Hele-Shaw flow}

\markboth{X. Feng and H. Wu}{Adaptive methods for the Cahn-Hilliard equation}

\author{
Xiaobing Feng\thanks{Department of Mathematics, The University
of Tennessee, Knoxville, TN 37996, U.S.A. ({\tt xfeng@math.utk.edu}).
The work of this author is partially supported by the NSF grant DMS-0410266.}
\and
Haijun Wu \thanks{Department of Mathematics, Nanjing University, Jiangsu,
210093, PR China. ({\em hjw@nju.edu.cn}). The work of this author is
partially supported by the China NSF grant 10401016, 
by the China National Basic Research Program grant 2005CB321701, and by 
the Natural Science Foundation of Jiangsu Province under the grant BK2006511.}
}

%\date{August 26, 2006}
\maketitle

\begin{abstract}
This paper develops a posteriori error estimates of residual type
for conforming and mixed finite element approximations of the fourth order
Cahn-Hilliard equation $u_t+\De\bigl(\eps \De u-\eps^{-1}
f(u)\bigr)=0$. It is shown that the {\it a posteriori} error bounds depends on
$\eps^{-1}$ only in some low polynomial order, instead of
exponential order. Using these a posteriori error estimates, we
construct an adaptive algorithm for computing the solution of the
Cahn-Hilliard equation and its sharp interface limit, the Hele-Shaw
flow. Numerical experiments are presented to show the robustness and
effectiveness of the new error estimators and the proposed adaptive
algorithm.
\end{abstract}



%%%%%%%%%%%%%%%%
\section{Introduction}\label{sec-1}
In this paper we derive a posteriori error estimates and
develop an adaptive algorithm based on the error estimates
for conforming and mixed finite element approximations of the
following Cahn-Hilliard equation and its sharp interface limit
known as the Hele-Shaw flow \cite{alikakos94,pego89}
\begin{eqnarray}
u_t+\De\bigl(\eps\De u-\frac{1}{\eps} f(u)\bigr)
&=& 0\quad\mbox{ in }\Om_T:=\Om\times (0,T),\label{e1.1}\\
\frac{\partial u}{\partial n} =\frac{\partial}{\partial n}
\bigl(\eps\De u-\frac{1}{\eps} f(u)\bigr) &=& 0\quad\mbox{ in
}\partial\Om_T
:=\partial\Om\times (0,T),\label{e1.2}\\
u &=& u_0\quad\mbox{in }\Om\times \{0\}, \label{e1.3}
\end{eqnarray}
where $\Om\subset \mathbf{R}^N$ $(N=2,3)$ is a bounded domain with
$C^{2}$ boundary $\partial\Om$  or a convex polygonal domain.
$T>0$ is a fixed constant, and $f$ is the derivative of a smooth
double equal well potential taking its global minimum value $0$ at
$u=\pm 1$. A well known example of $f$ is
\[
f(u):=F'(u)\quad\mbox{and}\quad F(u)=\frac{1}{4}(u^2-1)^2.
\]
For the notation brevity, we shall suppress the super-index $\eps$ on $u^\eps$
throughout this paper except in Section \ref{sec-5}.

The equation (\ref{e1.1}) was originally introduced by Cahn and Hilliard
\cite{cahn58} to describe the complicated phase separation and coarsening
phenomena in a melted alloy that is quenched to a temperature at which only
two different concentration phases can exist stably. The Cahn-Hilliard has
been widely accepted as a good (conservative) model to describe the phase
separation and coarsening phenomena in a melted alloy.
The function $u$ represents the concentration of one of
the two metallic components of the alloy.
The parameter $\varepsilon$ is an ``interaction length", which is
small compared to the characteristic dimensions on the laboratory
scale. Cahn-Hilliard equation (\ref{e1.1}) is a special case of a
more complicated phase field model for solidification of a pure
material \cite{caginalp86,fix83,langer86}. For the physical
background, derivation, and discussion of the Cahn-Hilliard equation
and related equations, we refer to \cite{allen79,alikakos94,bates93,cahn58,
chen96,elliott86,novick98,novick00} and the references therein. It should be
noted that the Cahn-Hilliard equation (\ref{e1.1}) can also be regarded as
the $H^{-1}$-gradient flow for the energy functional \cite{fife00}
\begin{equation}
{\mathcal J}_\varepsilon(u):=\int_\Om \Bigl[\,\frac12|\nabla
u|^2+\frac{1}{\eps^2} F(u)\,\Bigr]\, dx. \label{e4}
\end{equation}

In addition to its application in phase transition, the Cahn-Hilliard
equation (\ref{e1.1}) has also been extensively studied in the
past due to its connection to the following free boundary problem,
known as the Hele-Shaw problem and the Mullins-Sekerka problem
\begin{alignat}{2}
\De w&=0 &&\quad\mbox{in } \Om\setminus\G_t,\, t\in[0,T]\, ,
\label{eqn-1.6}\\
\frac{\p w}{\p n} &=0& &\quad\mbox{on }\p\Om,\, t\in[0,T]\, ,
\label{eqn-1.7}\\
w &=\sigma\kappa &&\quad\mbox{on }\G_t,\, t\in[0,T]\, ,
\label{eqn-1.8}\\
V&=\frac12\Bigl[\frac{\p w}{\p n}\Bigr]_{\G_t} &&\quad\mbox{on }\G_t,
\, t\in[0,T]\, , \label{eqn-1.9}\\
\G_{0} &=\G_{00} &&\quad\mbox{when } t=0\, . \label{eqn-1.10}
\end{alignat}
Here
\[
\sigma=\int_{-1}^1\sqrt{\frac{F(s)}{2}}\, {\rm d}s\, .
\]
$\kappa$ and $V$ are, respectively, the mean curvature and the
normal velocity of the interface $\G_t$, $n$ is the unit outward
normal to either $\p\Om$ or $\G_t$, $[\frac{\p w}{\p n}]_{\G_t}
:=\frac{\p w^+}{\p n}-\frac{\p w^-}{\p n}$, and $w^+$ and $w^-$
are respectively the restriction of $w$ in $\Om^+_t$ and
$\Om^-_t$, the exterior and interior of $\G_t$ in $\Om$.

Under certain assumption on the initial datum $u_0$,
it was first formally proved by Pego~\cite{pego89} that, as
$\varepsilon\searrow 0$, the function $w^\eps:=-\eps\De u^\eps
+\eps^{-1} f(u^\eps)$, known as the chemical potential, tends to
$w$, which, together with a free boundary $\G:=\cup_{0\leq t\leq
T}(\G_t\times\{t\})$ solves \eqref{eqn-1.6}-\eqref{eqn-1.10}. Also
$u^\eps\rightarrow\pm 1$ in $\Om_t^\pm$ for all $t\in [0,T]$, as
$\eps\searrow 0$. The rigorous justification of this limit was
carried out by Alikakos, Bates and Chen in \cite{alikakos94} under
the assumption that the above Hele-Shaw (Mullins-Sekerka) problem
has a classical solution. Later, Chen~\cite{chen96} formulated a
weak solution to the Hele-Shaw (Mullins-Sekerka) problem and showed,
using an energy method, that the solution of
(\ref{e1.1})-(\ref{e1.3}) approaches, as $\eps\searrow 0$, to a weak
solution of the Hele-Shaw (Mullins-Sekerka)
problem. One of a consequences of the connection between the Cahn-Hilliard
equation and the Hele-Shaw flow is that for small $\varepsilon$
the solution to (\ref{e1.1})-(\ref{e1.3}) equals $\pm 1$ in the
two bulk regions of $\Om$ which is separated by a thin layer
(called {\em diffuse interface}) of width $O(\epsilon)$. As
expected, the solution has a sharp moving front over the
transition layer.

Another motivation for developing efficient adaptive numerical methods for
the Cahn-Hilliard equation is its applications far beyond its original role in
phase transition. The Cahn-Hilliard equation is indeed a fundamental
equation and an essential building block in the {\em phase field theory}
for moving interface problems (cf. \cite{mcfadden02}), it is often combined
with other fundamental equations of mathematical physics such as the
Navier-Stokes equation (cf. \cite{feng05,jacqmin99,ls03} and the references
therein) to be used as diffuse interface models for describing various interface
dynamics, such as flow of two-phase fluids, from various applications.

The primary numerical challenge for solving the Cahn-Hilliard
equation results from the presence of the small parameter
$\varepsilon$ in the equation, so the equation is a {\em singular}
perturbation of the biharmonic heat equation. Numerically to resolve the
thin transition region of width $O(\varepsilon)$, one has to use very
fine meshes in the region. Considering the fact that away from the
transition region the solution equals $\pm 1$, it is natural
to use adaptive meshes, rather than uniform meshes, to compute
the solution. As far as the error analysis concerns, the main
difficulty is to derive a priori and {\em a posteriori}
error estimates which depends on $\frac{1}{\varepsilon}$ only in
(low) polynomial order, rather than exponential order which is
the case if the standard Gronwall's inequality type argument
is used to derive the error estimates \cite{barrett01,du91,elliottc89,
elliott89}. Recently, Feng and Prohl \cite{XA2,XA3,XA4} were able to
overcome this difficulty and established polynomial order a priori
error estimates for mixed finite element approximations
of the Cahn-Hilliard equation and related phase field equations. Based on these
new error estimates, they then proved convergence of the numerical solutions
of the phase field equations to the solutions of their respective sharp
interface limits as mesh sizes and the parameter $\eps$ all tend to zero.
The main idea of \cite{XA2,XA3} is to use a spectral estimate result
of Alikakos and Fusco \cite{alikakos93} and Chen \cite{chen94}
for the linearized Cahn-Hilliard operator to handle the nonlinear term
in the error equation.  Very recently, this idea was also
used by Kessler, Nochetto and Schmidt \cite{nochetto03} and by Feng and Wu
\cite{feng_wu03} to obtain {\em a posteriori} error estimates, which
depend on $\frac{1}{\eps}$ in some low polynomial order, for
finite element approximations of the Allen-Cahn equation.

The goal of this paper is
to develop {\em a posteriori} error estimates for
conforming and mixed finite element approximations
of the Cahn-Hilliard equation in the spirit of
\cite{feng_wu03}. First, using the idea of continuous dependence 
we derive some residual type a posteriori error
estimates, which depend on $\frac{1}{\eps}$ only in low polynomial
orders, for the conforming finite element approximations
and the mixed finite element approximations.
To avoid many technicalities and to present the idea, we only consider
semi-discrete (in spatial variable) approximations in this paper. For the
time discretization, we appeal to the stiff ODE solver NDF \cite{sr97}
which is a modification of BDF for temporal integration. Then, using the
a posteriori estimates as error indicators we propose an
adaptive algorithm for approximating the Cahn-Hilliard equation
and its sharp interface limit, the Hele-Shaw flow.
As in \cite{feng_wu03}, the technique and analysis of this paper
for deriving a posteriori error estimates are problem-independent
and method-independent, hence, they are applicable to
a large class of evolution problems and their numerical approximations
obtained by any (numerical) discretization method including finite difference,
finite element, finite volume and spectral methods. We also remark that
the adaptive finite element algorithm of this paper is based on
the method of lines approach, we refer to \cite{af88, bfs01, ej95} and
the references therein for a detailed exposition on the approach
for other types of problems, and to \cite{ej95,verfurth98}
and the references therein for a detailed discussions about
adaptive algorithms based on other approaches such as
discontinuous Galerkin methods and space-time finite element methods.

The paper is organized as follows: In Section \ref{sec-2} we establish
continuous dependence estimates for the Cahn-Hilliard equation in both
standard and mixed formulations, and present some abstract frameworks
for deriving a posteriori error estimates based on the idea of
continuous dependence. In Section \ref{sec-3} we derive some a posteriori error
estimates for conforming finite element approximations
and for the Ciarlet-Raviart mixed finite element approximations of
the Cahn-Hilliard
equation using the continuous dependence estimates and the abstract
frameworks of Section \ref{sec-2}. In Section \ref{sec-4} we
propose an adaptive finite element algorithm using the a posteriori
error estimates of Section \ref{sec-3} as error indicators for
refining or coarsening the mesh. In Section \ref{sec-5} we
establish some a posteriori error estimates for using the
conforming and mixed finite element methods to approximate the Hele-Shaw flow.
Finally, in Section \ref{sec-6} we present several numerical tests
to show the robustness and effectiveness of the
proposed error estimators and the adaptive algorithm.

\section{Continuous dependence and a posteriori error estimates} \label{sec-2}

In this section, we first establish some continuous dependence (on
nonhomogeneous force term and on initial condition) estimates for the
Cahn-Hilliard problem \reff{e1.1}-\reff{e1.3} in both standard and mixed
formulations. We then present an abstract framework
for deriving a posteriori error estimates for mixed numerical
approximations of general evolution equations. Our goal is
to derive a posteriori error estimates which depend on $\frac{1}{\eps}$
only in some low polynomial order. It is easy to show that (cf. Section
\ref{sec-2.1} ) if one uses the standard perturbation and Gronwall's
inequality techniques to derive a priori or a posteriori error estimates,
the error bounds will depend on $\frac{1}{\eps}$ exponentially, hence,
such estimates are not useful for small $\varepsilon$. To overcome the
difficulty, we appeal to a spectrum estimate result, due to Alikakos and
Fusco \cite{alikakos93} and Chen \cite{chen94},
for the linearized Cahn-Hilliard operator, and prove a continuous
dependence estimate, which depends on $\frac{1}{\eps}$
in some low polynomial order, for the  Cahn-Hilliard equation.
Such a continuous dependence estimate is the key for us to establish the
desired a posteriori error estimates in the next section.

Throughout this paper, the standard space, norm and inner product
notation are adopted. Their definitions can be found in
\cite{brenner94,ciarlet78}. In particular, $(\cdot,\cdot)$ denotes
the standard $L^2$-inner product, and $H^k(\Om)$ stands
for the usual Sobolev spaces. Also, $C$ are used to denote a generic
positive constant which is independent of $\varepsilon$ and the
mesh sizes.

\subsection{Continuous dependence estimates} \label{sec-2.1}

Introduce the space
\[H^2_E(\Om)=\set{\psi\in H^2(\Om);\;\frac{\partial\psi}{\partial n}=0
\text{ on } \partial\Om}.
\]
 We recall that the variational formulation
of \eqref{e1.1}--\eqref{e1.3} is defined by seeking $u\in
H^2_E(\Om)$ such that
\begin{align}\label{eu}
\langle u_t, \psi\rangle +\eps\bigl(\De u, \De \psi\bigr)
+\frac{1}{\eps} \bigl(\na(f(u)),\na \psi\bigr)
&=0 \quad\forall \psi\in  H^{2}(\Om),\;t\in [0,T], \\
u(0) &=u_0\in H^2_E(\Om) .  \label{eu1}
\end{align}
It is proved in \cite{elliottc89} that such a solution $u$ exists and
\[u\in L^\infty((0,T);
H^2_E(\Om))\cap L^2((0,T);H^4(\Om))\cap H^1((0,T);L^2(\Om)).\]
 For physical reason, unless mentioned
otherwise, we assume that $|u_0|\leq 1$ in this paper.

Let $v(t)\in H^2_E(\Om)$ be a perturbation of $u$ satisfying
\begin{alignat}{2}\label{ev}
   \dual{v_t}{\psi} + \eps\bigl(\De v, \De \psi\bigr)
+\frac{1}{\eps} \bigl(\na(f(v)),\na \psi\bigr)
= \dual{r(t)}{\psi} &&&\quad\forall \psi\in  H^2_E(\Om),\;t\in [0,T], \\
   v(0) =v_0\in H^2_E(\Om), &&& \label{ev1}
\end{alignat}
where $r(t)\in\shmlt:=(H^2_E(\Om))^*$ (the dual space of
$H^2_E(\Om)$) is the residual of $v(t)$, i.e., the perturbation of
the right-hand side of \eqref{e1.1}. $\dual{\cdot}{\cdot}$ denotes
the dual product on $\shmlt\times H^2(\Om)$. We assume that
$\dual{r(t)}{1} = 0$, and define
\begin{equation}\label{em1norm}
\norm{r(t)}{\shmt} = \sup_{0\neq\psi \in
H^2_E(\Om)}\frac{\dual{r(t)}{\psi}}{\norm{\psi}{H^2}} .
\end{equation}
Let $L^2_0(\Om)=\set{\psi\in L^2(\Om);\;\int_\Om\psi dx = 0}$. Define
$\De^{-1}: L^2_0(\Om)\rightarrow H^1(\Om) \cap L^2_0(\Om)$ to be the
inverse of the Laplacian $\De$, that is, for any $\psi\in
L^2_0(\Om)$, $\De^{-1} \psi\in H^1(\Om) \cap L^2_0(\Om)$ is defined
by
\[
\bigl(\na(\De^{-1}\psi),\na\eta\bigr) =-(\psi,\eta)\qquad\forall
\eta\in H^1(\Om).
\]
From the standard regularity theory of elliptic problems, one concludes
that $\De^{-1} \psi\in H^2_E(\Om)$ and
\begin{equation}\label{ereg}
    \norm{\De^{-1} \psi}{H^2(\Om)}\le
    C\normL{\psi}.
\end{equation}
Let $w(t):=v(t) - u(t)$. We also assume that $w(0)=v_0-u_0\in
L^2_0(\Om)$. Then, from $\int_\Om w(t) dx=\int_\Om w(0) dx$, it is clear
that $w(t)\in L^2_0(\Om)$. Subtracting equation \eqref{eu} from
equation \eqref{ev} gives
\begin{equation}\label{eerr}
   \dual{w_t}{\psi} + \eps\bigl(\De w, \De \psi\bigr)
+\frac{1}{\eps} \bigl(\na(f(v)-f(u)),\na \psi\bigr) =
\dual{r(t)}{\psi} \quad\forall \psi\in  H^2_E(\Om).
\end{equation}

Next, we give two estimates on $u-v$ in terms of $r$ and $u_0-v_0$
for the Cahn-Hilliard equation. The first estimate holds without
any constraint on either the initial condition or the residual of the
perturbation problem, but the estimate depends on
$\frac{1}{\varepsilon}$ { \em exponentially}. The second one,
which depends on $\frac{1}{\varepsilon}$ only in a low polynomial
order, holds provided that the perturbations of the initial condition and the
right-hand side are small.

\begin{proposition}\label{prop2.1}
Let $u$ and $v$ be the weak solutions of \eqref{eu}-\eqref{eu1}
and \eqref{ev}-\eqref{ev1}, respectively. Then it holds that for
$t\in [0, T]$
\begin{equation}\label{e2.2}
\begin{split}
&\normL{\na \De^{-1} (v(t) - u(t))}^2 +\eps\int_0^t
\exp\Bigl(\frac{4(t-s)}{\eps^3}\Bigr)\normL{\na (v(s)-u(s))}^2\,ds \\
&\hskip 0.3in
\leq \exp\Bigl(\frac{4t}{\eps^3}\Bigr)\, \normL{\na \De^{-1}
(v_0-u_0)}^2 +\frac{C}{\eps}\int_0^t
\exp\Bigl(\frac{4(t-s)}{\eps^3}\Bigr)\norm{r(s)}{\shmt}^2\, ds .
\end{split}
\end{equation}

\end{proposition}

\begin{proof}
Setting $\psi =-\De^{-1} w$ in \eqref{eerr} we get
\begin{equation}\label{e2.3}
\frac12\frac{d}{dt}\normL{\na \De^{-1}w}^2 + \eps\normL{\na w}^2
+\frac1\eps\bigl(f(v)-f(u), w \bigr) = -\dual{r}{\De^{-1}w}.
\end{equation}
From the definition of $\De^{-1}$ it follows
\begin{equation}\label{ewl2}
    \normL{w}^2=\bigl(\na(\De^{-1}w),\na w\bigr)\le
    \normL{\na(\De^{-1}w)}\normL{\na
    w}.
\end{equation}
Hence,
\begin{align*}
 \frac1\eps\bigl(f(v)-f(u), w\bigr)&=\frac1\eps\bigl(f^\prime(\xi)w, w\bigr)
 =\frac1\eps\bigl((3\xi^2-1)w, w\bigr) \ge -\frac1\eps\normL{w}^2\\
 &\ge -\frac{\eps}{4}\normL{\na w}^2 -
 \frac1{\eps^3}\normL{\na(\De^{-1}w)}^2 .
\end{align*}
Similarly,
\begin{align*}
    -\dual{r}{\De^{-1}w}&\le\norm{r}{\shmt}\norm{\De^{-1}w}{H^2}
    \le C\norm{r}{\shmt}\normL{w}
    \le C\eps\norm{r}{\shmt}^2 + \frac1\eps\normL{w}^2\\
    &\le C\eps\norm{r}{\shmt}^2+\frac{\eps}{4}\normL{\na w}^2 +
 \frac1{\eps^3}\normL{\na(\De^{-1}w)}^2 .
\end{align*}
Combining the above two estimates and \eqref{e2.3} we obtain
\begin{equation*}
\frac{d}{dt}\normL{\na \De^{-1}w}^2 + \eps\normL{\na w}^2
\le\frac4{\eps^3}\normL{\na(\De^{-1}w)}^2 + C\eps\norm{r}{\shmt}^2 .
\end{equation*}
Finally, the desired estimate \eqref{e2.2} follows from an application
of the Gronwall's inequality. The proof is complete.
\end{proof}

\begin{remark}
Clearly, the above continuous dependence estimates are only useful when
$t=O(\eps^3)$. However, the estimate is sharp {\em if} no assumptions
on the solutions $u$ and $v$ are assumed because the Cahn-Hilliard
equation does exhibit a fast initial transient regime for times
of order $O(\varepsilon^3)$, until interfaces develop \cite{cahn58,alikakos94}.
\end{remark}

To improve estimates \eqref{e2.2}, we need to confine ourself to consider
solutions $u$ and $v$ which have certain profiles. Specifically, we need
the helps of the following
three lemmas. The first lemma gives an a priori estimate for
solutions of a Bernoulli type nonlinear ordinary differential
inequality. Its proof can be found in \cite{feng_wu03}.

\begin{lemma}\label{lem2.1} Suppose that $n>1$, $y(t)$ and $\lambda(t)$ are
nonnegative functions satisfying
\begin{equation} \label{e2.4}
y^\prime(t)\leq \lambda(t)\, (y(t))^n + a(t) y(t) + b(t) \qquad
\forall t\in [0,T] \, .
\end{equation}
Define $\rho(t)=\int_0^te^{-\int_0^s a(\tau)\,d\tau} b(s)\,ds$ and
$\bar{\rho}(t)=\max_{0\leq s\leq t}\rho(s)$, then there holds for
$t\in [0,T^*)$
\begin{equation}\label{e2.6}
y(t) \leq \frac{[y(0) +\bar{\rho}(t)]\, e^{\int_0^t
a(s)\,ds}}{\zeta(t)^{\frac{1}{n-1}}}+[\rho(t)-\bar{\rho}(t)]e^{\int_0^t
a(s)\,ds} ,
\end{equation}
where
\begin{equation*}
\zeta(t) = 1-(n-1)\,[y(0)+\bar{\rho}(t)]^{n-1} \int_0^t \,
\lambda(s)\, e^{(n-1)\int_0^s a(\tau)\, d\tau}\, ds,
\end{equation*}
and $T^*$ is the largest positive number in $[0,T]$ such that
$\zeta(t)\geq 0$ .
\end{lemma}

The second lemma cites a spectrum estimate result of Alikakos and Fusco
\cite{alikakos93} and Chen \cite{chen94} for the following
linearized Cahn-Hilliard operator at the solution of \eqref{e1.1}-\eqref{e1.3}
\begin{equation}\label{Ch}
{\mathcal L}_{CH}:=\De\bigl(\eps \De -\frac{1}{\eps} f^\prime(u)
I\bigr),
\end{equation}
where $I$ stands for the identity operator.

\begin{lemma}\label{lem2.2}
Let $\lambda_{CH}$ denote the smallest eigenvalue of ${\mathcal
L}_{CH}$, assume that the solution $u$ satisfies the $\tanh$ profile
described in \cite{chen94} (cf. (1.10) on page 1374 and Theorem 1.1
on page 1375 of \cite{chen94}). Then
there exists $0<\varepsilon_0<1$ and an
$\varepsilon$-independent positive constant $C_0$ such that
$\lambda_{CH}$ satisfies
\begin{equation*}
\lambda_{CH} \equiv \inf_{0\not\equiv\psi \in H^{1}(\Om)\cap
L^2_0(\Om)} \frac{\varepsilon\, \normL{\nabla \psi}^2 +
\frac{1}{\varepsilon}\, (f'(u) \psi, \psi)}{ \normL{\na \De^{-1}
\psi}^2} \geq -C_0 \qquad \forall \varepsilon\in (0, \varepsilon_0].
\end{equation*}
\end{lemma}

\begin{remark} Since the proof of the above estimate is based on
the convergence result of \cite{caginalp98}, which says that the
solution of the Cahn-Hilliard problem (\ref{e1.1})-(\ref{e1.3})
for certain class of initial conditions converges to the classical
solution of the free boundary problem (\ref{eqn-1.6})-(\ref{eqn-1.10})
as $\varepsilon\rightarrow 0$, hence, the proof suggests that the
validity of the above estimate also depends on the choice of
the initial conditions. As far as we know it is an open question
whether the estimate still holds for ``general" initial data
(see Remark 2.3 of \cite{caginalp98} for more discussions). This is the
reason why the subsequent a posteriori error estimates of this paper
are established under this initial condition constraint.
\end{remark}

The third lemma gives an estimate which are useful for the
subsequent analysis.
\begin{lemma}\label{lem2.3}
Let $0<\de<2$, then there exits a positive constant $C$ which is independent
of $\varepsilon$ and $\de$ such that for any $w\in H^1(\Om)\cap L^2_0(\Om)$
there holds
\begin{equation}\label{ew3}
 \frac{1}\eps\int_\Om \abs{w}^3\, dx \le
 \frac1{2\eps}\norm{w}{L^4}^4 + \frac{\eps^4}4\normL{\na w}^2 +
C\de\eps^{4-\frac{20}\de}\normL{\na\De^{-1}w}^{\frac{16+2(N-2)\de}{(2+N)\de}}.
\end{equation}
\end{lemma}
\begin{proof}
Recall the Young's inequality
\[ab\le \frac{q-1}{q}a^{\frac{q}{q-1}}+\frac{b^q}{q},\qquad a, b>0, q>1.\]
Hence,
\begin{equation}\label{eyoung}
    ab\le
    a^{\frac{q}{q-1}}+\big(1-\frac1q\big)^q\frac{b^q}{q-1}
    \le a^{\frac{q}{q-1}}+e^{-1}\frac{b^q}{q-1} .
\end{equation}
Then for $2< p<3$
\begin{align*}
    \abs{w}^3=\left(\frac{\abs{w}^4}{2}\right)^{\frac{3-p}{4-p}}
    2^{\frac{3-p}{4-p}}\abs{w}^{\frac{p}{4-p}}\le\frac{\abs{w}^4}{2}+C\abs{w}^p,
\end{align*}
therefore,
\begin{equation}\label{ew3a}
   \frac{1}\eps \int_\Om |w|^3\, dx\le \frac1{2\eps}\norm{w}{L^4}^4
   +\frac{C}{\eps}\norm{w}{L^p}^p .
\end{equation}
Since $w\in  H^1(\Om)\cap L^2_0(\Om)$, it follows from the Sobolev
inequality and \eqref{ewl2} that
\begin{align*}
\norm{w}{L^p} %\le C\norm{w}{H^{\frac{N(p-2)}{2p}}}
\le\normL{w}^{1-\frac{N(p-2)}{2p}}\normL{\na w}^{\frac{N(p-2)}{2p}}
\le C\normL{\na\De^{-1}w}^{\frac{2p-N(p-2)}{4p}}\normL{\na
w}^{\frac{2p+N(p-2)}{4p}}.
\end{align*}
Let $p=\frac{8+2N-2\de}{2+N}=2+\frac{2(2-\de)}{2+N}$, we have
\begin{align*}
\frac{1}\eps\norm{w}{L^p}^p&
\le C\eps^{-5+\de} \left(\frac{\eps^4}{4}
\normL{\na w}^2\right)^{\frac{4-\de}{4}}
\normL{\na\De^{-1}w}^{\frac{8+(N-2)\de}{2(2+N)}}.
\end{align*}
From inequality \eqref{eyoung} with $q=\frac4{\de}$ we obtain
\begin{align*}
\frac{1}\eps\norm{w}{L^p}^p&\le \frac{\eps^4}{4}\normL{\na w}^2
+C\de\eps^{4-\frac{20}{\de}}
\normL{\na\De^{-1}w}^{\frac{16+2(N-2)\de}{(2+N)\de}}.
\end{align*}
\eqref{ew3} now follows from combining the above estimate and \eqref{ew3a}.
The proof is complete.
\end{proof}

We are now ready to state our first main result of this section.

\begin{proposition}\label{prop2.3}
Suppose that $|u_0|, |v_0|\leq 1$, $\eps_0$ and $C_0$ be the same
as in Lemma \ref{lem2.2}. Let $u$ and $v$ be the solutions of
\eqref{eu}-\eqref{eu1} and \eqref{ev}-\eqref{ev1}, respectively.
Then, for any $\varepsilon\in (0,\varepsilon_0]$, there exists a
positive constant $C$, which is independent of $\eps$ and $t$,
such that there holds
\begin{equation}\label{e2.7}
\begin{split}
&\normL{\na\De^{-1}(v(t)-u(t))}^2\\
&\qquad+\int_0^t\,\Bigl(\eps^4\normL{\nabla(v(s)-u(s))}^2
+\frac{1}{\eps}\norm{v(s)-u(s)}{L^4}^4\Bigr)e^{(2C_0+8)(t-s)}\, ds\\
& \leq \frac{1}{\xi(t)} \normL{\na\De^{-1} (v_0-u_0)}^2e^{(2C_0+8)t}\\
&\qquad+\Bigl[1+\frac{1}{\xi(t)}\Bigr]C \eps^{-2}\int_0^t\,
\norm{r(s)}{\shmt}^2 e^{(2C_0+8)(t-s)}\, ds
\end{split} \end{equation}
for all $t\in [0,T^*)$. Here
\begin{equation}\label{e2.8a}
\begin{split}
\xi(t):=&1- C \eps^{-\frac{5(2+N)}2} e^{(2C_0+8)t}\times\\
&\qquad\left\{\normL{\na \De^{-1} (v_0-u_0)}^2
+\eps^{-2}\int_0^t\, \norm{r(s)}{\shmt}^2\, e^{-(2C_0+8)s}\, ds\right\} ,
\end{split}
\end{equation}
and $T^*\in [0,T]$ satisfying $\xi(T^*)> 0$.
\end{proposition}

\begin{proof}
Let $w:=v-u$, from \eqref{e2.3} and the identities
\begin{align}\label{efvu}
f(v) - f(u) &= f^\prime(u) w + w^3 + 3 u w^2,\\
\bigl( f(v) - f(u), w \bigr) &=\int_\Om f^\prime(u)\, w^2\, dx +
\norm{w}{L^4}^4 + 3\int_\Om u\, w^3\, dx,  \nonumber
\end{align}
and the fact that $\norm{u}{L^\infty}\leq C$ (cf. \cite{caffarelli95,XA3})
we have
\begin{equation}\label{e2.9}
\begin{split}
&\frac12 \frac{d}{d t} \normL{\na\De^{-1} w}^2 + \frac{1}{\eps} \norm{w}{L^4}^4
+\eps\normL{\nabla w}^2 +\frac{1}{\eps} \int_\Om f^\prime(u)\, w^2\, dx  \\
&\quad =-\frac{3}{\eps} \int_\Om u\,w^3\, dx - \dual{r}{\De^{-1} w}
\le\frac{C}{\eps} \int_\Om |w|^3\, dx + \frac{C}{\eps^2}\norm{r}{\shmt}^2
+\eps^{2}\normL{w}^2.
\end{split}
\end{equation}

To bound the fourth term on the left-hand side of \eqref{e2.9}
from below, we employ the spectrum estimate of Lemma \ref{lem2.2}.
In order to keep a portion of $\normL{\nabla w}^2$ on the
left-hand side, we apply the spectrum estimate with a scaling
factor $(1-\eps^{3})$.
\begin{equation*}
\begin{split}
&\eps\normL{\nabla w}^2 +\frac{1}{\eps}
\int_\Om f^\prime(u)\, w^2 dx -\eps^{2}\normL{w}^2\\
&\hskip 0.2in
=\eps^{3}\, \Bigl[\, \eps \normL{\nabla w}^2
+\frac{1}{\eps} \int_\Om \big(3u^2-2\big)\, w^2 dx \Bigr]
+ (1-\eps^{3})\, \Bigl[\, \eps \normL{\nabla w}^2
+\frac{1}{\eps} \bigl(f^\prime(u)w, w \bigr)\,\Bigr]  \\
&\hskip 0.2in
\geq \eps^4 \normL{\nabla w}^2 - C_0\normL{\na\De^{-1} w}^2
-2\eps^{2} \normL{w}^2 .
\end{split}
\end{equation*}
Since
\begin{align*}
2\eps^2\normL{w}^2\le 2\eps^2\normL{\na w}\normL{\na\De^{-1}w}\le
\frac{\eps^4}4\normL{\na w}^2 + 4\normL{\na\De^{-1}w}^2,
\end{align*}
we have
\begin{equation}\label{e2.11}
\eps \normL{\nabla w}^2 +\frac{1}{\eps}
\int_\Om f^\prime(u)\, w^2 dx -\eps^{2}\normL{w}^2\\
\geq \frac{3\eps^4}4 \normL{\nabla w}^2 -(C_0+4)\normL{\na\De^{-1} w}^2 .
\end{equation}
Combining \eqref{e2.11}, \eqref{ew3}, and \eqref{e2.9} we obtain
\begin{equation}\label{e2.12}
\begin{split}
\frac{d}{d t}\normL{\na \De^{-1} w}^2
&\le C\de\eps^{4-\frac{20}{\de}}
\normL{\na\De^{-1}w}^{\frac{16+2(N-2)\de}{(2+N)\de}}
+(2C_0+8)\normL{\na\De^{-1} w}^2 \\
&\hskip .5in
+C\eps^{-2}\norm{r}{\shmt}^2  -\eps^4\normL{\nabla w}^2
-\frac{1}{\eps} \norm{w}{L^4}^4 ,
\end{split}
\end{equation}
where $0<\de<2$.

Now, set
\begin{gather*}
y(t):=\normL{\na \De^{-1} w}^2,\quad a:=2C_0+8,\quad
\lambda:=C\de\eps^{4-20/\de} ,
\quad n:=\frac{8+(N-2)\de}{(2+N)\de}, \\
b(t):= C\eps^{-2}\norm{r}{\shmt}^2 -\eps^4\normL{\nabla w}^2
-\frac{1}{\eps} \norm{w}{L^4}^4, \quad
 \rho(t):=\int_0^t e^{-(2C_0+8)s}b(s)\,ds,
\end{gather*}
then
\[
0\leq\bar{\rho}(t)\leq C\int_0^t
e^{-(2C_0+8)s}\eps^{-2}\norm{r(s)}{\shmt}^2\, ds\, .
\]
It follows from Lemma \ref{lem2.1} that there exists $T^*\in
(0,T]$ such that
\begin{equation}\label{ey} y(t)\leq
\frac{\bigl(y(0)+\bar{\rho}(t)\bigr)\, e^{(2C_0+8)t} } { (
\zeta(t))^{\frac{1}{n-1}} } + e^{(2C_0+8)t} \rho(t)
\end{equation}
for all $t\in (0,T^*)$, where
\[\zeta(t)=1- \frac{\lambda}{2C_0+8}\,
\bigl[y(0)+\bar{\rho}(t)\bigr]^{n-1} \bigl[\, e^{(2C_0+8)(n-1)t}
-1 \,\bigr]\,.\]
Moreover, since
\begin{align*}
(\zeta(t))^{\frac{1}{n-1}}&\ge \left(1- \frac{\lambda}{2C_0+8}\,
\bigl[y(0)+\bar{\rho}(t)\bigr]^{n-1}
\, e^{(2C_0+8)(n-1)t}\right)^{\frac1{n-1}}\\
&\ge 1- \left(\frac{\lambda}{2C_0+8}\right)^{\frac1{n-1}}\,
\bigl[y(0)+\bar{\rho}(t)\bigr] \, e^{(2C_0+8)t} ,
\end{align*}
then there exists a positive constant $C$ independent of $\eps$ and $\de$
such that
\begin{align*}
(\zeta(t))^{\frac{1}{n-1}}\ge 1- C\eps^{-\frac{(2+N)(5-\de)}{(2-\de)}}\,
\bigl[y(0)+\bar{\rho}(t)\bigr] \, e^{(2C_0+8)t}.
\end{align*}
The estimate \eqref{e2.7} now follows from combining the above inequality
and \eqref{ey} and letting $\de\to 0$. The proof is complete.
\end{proof}

\begin{remark}
In the above proof we have used the boundedness property of the solution of
the Cahn-Hilliard problem \eqref{e1.1}--\eqref{e1.3}, which will be used
a couple more times later in the paper. The references we cited for the
property are \cite{caffarelli95,XA3}. However, we like to point out that
the assertion was proved in \cite{caffarelli95} under the assumption that
the derivative $f(u)=F'(u)$ of the potential $F$ is linear outside
a bounded interval, which is not the case for the potential
$F(u)=\frac14(u^2-1)^2$ used in this paper. Although we believe
the boundedness of the solution in the case of the above potential
also holds, we have not found a (direct) proof in the literature. On the
other hand, an indirect proof was given in \cite{XA3}
(see Lemma 2.2 of \cite{XA3}),
which uses the fact that the solution of the Cahn-Hilliard
problem (\ref{e1.1})-(\ref{e1.3}) converges to the classical
solution of the free boundary problem (\ref{eqn-1.6})-(\ref{eqn-1.10})
as $\varepsilon\rightarrow 0$.  As a result, the proof depends
on the choice of the initial conditions. Hence, as pointed out
in Remark 2.2, the subsequent a posteriori error estimates of this paper
are established under this initial condition constraint.
\end{remark}

\medskip
In order to assure the continuous dependence estimate of Proposition
\ref{prop2.3} hold on the whole interval $(0,T)$, we need to
impose a {\em smallness} constraint on the perturbations of the
initial condition and the right-hand side as described in the
following corollary.

\begin{corollary}\label{cor2.1}
Under the assumptions of Proposition \ref{prop2.3}, estimate
\eqref{e2.7} holds for $T^*=T$ if $v_0$ and $r$ satisfy the
following constraint
\begin{equation}\label{e2.13}
\begin{split}
&\left\{\normL{\na \De^{-1} (v_0-u_0)}^2 +\eps^{-2}\int_0^T\,
\norm{r(s)}{\shmt}^2\, e^{-(2C_0+8)s}\,
ds\right\}^{\frac12} \\
&\hskip 1in
\leq C^{-1} e^{-(C_0+4)T} \eps^{\frac{5(2+N)}{4}}
=\left\{
\begin{array}{ll}
        O(\eps^{5}) &\quad\mbox{if } N=2,\\
        O(\eps^{6.25}) &\quad\mbox{if } N=3.
        \end{array} \right.
\end{split}
\end{equation}
\end{corollary}

\begin{proof}
The assertion follows immediately from the fact that $\xi(T)> 0$
when \eqref{e2.13} holds.
\end{proof}

\begin{proposition}\label{prop2.2}
Under the assumptions of Corollary \ref{cor2.1}, there exists a
constant $C$ independent of $\eps$ such that for $t\in [0, T]$
\begin{equation}\label{e2.2.1}
\begin{split}
&\normL{v(t) - u(t)}^2 +\int_0^t\Big(\eps\normL{\De
(v(s)-u(s))}^2\\
&\hskip 1.5in+\frac{1}{\eps}\normL{(v(s)-u(s))\na
(v(s)-u(s))}^2\Big)\,ds \\
&
\leq \normL{v_0-u_0}^2 +\frac{C}{\eps^5\xi(t)}
\normL{\na \De^{-1} (v_0 - u_0)}^2e^{(2C_0+8)t} \\
&\hskip 0.9in
+\frac{C}{\eps^7}\left[1+\frac{1}{\xi(t)} \right]\int_0^t\,
\norm{r(s)}{\shmt}^2 e^{(2C_0+8)(t-s)}\, ds\,.
\end{split}
\end{equation}
\end{proposition}

\begin{proof}
Setting $\psi = w:=v(t) - u(t)$ in \eqref{eerr} gives
\begin{equation}\label{e2.3.1}
\frac12\frac{d}{dt}\normL{w}^2 + \eps\normL{\De w}^2
+\frac{1}{\eps}\bigl(\na(f(v)-f(u)),\na w \bigr) = \dual{r}{w}.
\end{equation}
From \eqref{efvu}, \eqref{ewl2}, and the fact that $\norm{u}{L^\infty}<C$
(cf. \cite{caffarelli95,XA3}) we get
\begin{align*}
\frac{1}{\eps}\bigl(&\na(f(v)-f(u)),\na w
\bigr)=\frac{1}{\eps}\bigl(\na(w^3+f^\prime(u) w +3u w^2),\na w \bigr)\\
&\qquad
=\bigl(3w^2\na w,\na w \bigr)
-\frac{1}{\eps}\bigl(f^\prime(u) w +3u w^2, \De w\bigr)\\
&\qquad
\ge \frac3\eps\normL{w\na w}^2-\frac{\eps}4\normL{\De w}^2
-\frac{C}{\eps^3}\Big(\norm{w}{L^2}^2+\norm{w}{L^4}^4\Big)\\
&\qquad
\ge \frac3\eps\normL{w\na w}^2-\frac{\eps}4\normL{\De w}^2
-\frac{C}{\eps^3}\Big(\frac1{\eps^2}\normL{\na\De^{-1}w}^2
+\eps^2\normL{\na w}^2+\norm{w}{L^4}^4\Big) .
\end{align*}
Combining this estimate and \eqref{e2.3.1} yields
\begin{align*}
&\frac12\frac{d}{dt}\normL{w}^2 + \frac{3\eps}4\normL{\De w}^2
+\frac{3}{\eps}\normL{w\na w}^2 \\
&\qquad
\le\frac{C}{\eps^5}\big(\normL{\na\De^{-1}w}^2
+\eps^4\normL{\na w}^2+\eps^2\norm{w}{L^4}^4\big)
+C\norm{r(s)}{\shmt}\normL{\De w}\\
&\qquad
\le\frac{C}{\eps^5}\big(\normL{\na\De^{-1}w}^2
+\eps^4\normL{\na w}^2+\eps^2\norm{w}{L^4}^4\big)
+\frac{C}{\eps}\norm{r(s)}{\shmt}^2+\frac{\eps}{4}\normL{\De w}^2.
\end{align*}
Here we have used the inequality $\norm{w}{H^2}=\norm{\De^{-1}\De w}{H^2}
\le C\normL{\De w}$ (cf. \eqref{ereg}) to derive the first inequality.
Therefore
\begin{align*}
\frac{d}{dt}\normL{w}^2 &+ \eps\normL{\De w}^2
+\frac{1}{\eps}\normL{w\na w}^2 \\
&\quad
\le \frac{C}{\eps^5}\big(\normL{\na\De^{-1}w}^2 +\eps^4\normL{\na
w}^2+\eps^2\norm{w}{L^4}^4\big) +C\eps^{-1}\norm{r(s)}{\shmt}^2.
\end{align*}
Integrating the above inequality over $[0, t]$ and using Proposition
\ref{prop2.3} and Corollary \ref{cor2.1} give \eqref{e2.2.1}.
The proof is complete.
\end{proof}

\subsection{Continuous dependence estimates for the mixed formulation}
\label{sec-2.2}

In this subsection we derive a continuous dependence estimate which is
analogous to \eqref{e2.7} for a {\em mixed} formulation of the
Cahn-Hilliard equation. It is well known that although at the differential
level the mixed weak formulation and the standard weak formulation are
equivalent, they are usually very different at the discrete level, i.e.,
the approximate
solutions obtained using these two variational formulations are
quite different. Indeed, it will be seen from the following estimate
that the mixed weak formulation results in two residual terms while
the standard weak formulation only gives one residual term, and in
general the combined effect of the former are not same as the effect
of the later.

Recall that \cite{XA3} the mixed formulation of problem
\eqref{eu}-\eqref{eu1} is defined by seeking a pair of functions
$(u(t),\vphi(t))\in [H^1(\Om)]^2$ such that
\begin{align}\label{e2.14}
\bigl(u_t,\psi\bigr)+ \bigl(\na \vphi,\na\psi\bigr) &=0
\quad\forall\psi\in H^1(\Om),\, t\in [0,T],\\
\eps \bigl(\na u,\na\chi \bigr)+\frac{1}{\eps} \bigl(f(u),
\chi\bigr) -\bigl(\vphi,\chi\bigr) &=0\quad\forall \chi\in
H^1(\Om),\, t\in [0,T],
\label{e2.15}\\
u(0) &= u_0\quad\mbox{in } \Om.  \label{e2.16}
\end{align}

We now consider a perturbation $(v(t),\phi(t))\in [H^1(\Om)]^2$ of
$(u(t),\vphi(t))$ defined by
\begin{align}\label{e2.17}
\bigl(v_t,\psi\bigr)+ \bigl(\na \phi,\na\psi\bigr)
&=\dual{r_1}{\psi}
\quad\forall\psi\in H^1(\Om),\, t\in [0,T],\\
\eps \bigl(\na v,\na\chi \bigr)+\frac{1}{\eps} \bigl(f(v),
\chi\bigr) -\bigl(\phi,\chi\bigr) &=\dual{\eps r_2}{\chi}
\quad\forall \chi\in H^1(\Om),\, t\in [0,T], \label{e2.18}\\
v(0) &= v_0\quad\mbox{in } \Om \label{e2.19}
\end{align}
for given ``residuals" $(r_1(t),r_2(t))\in [(H^1(\Om))^*]^2$ which
satisfy $\dual{r_1}{1}=\dual{r_2}{1}=0$. Introduce the following
norms of $r_j, j=1, 2$
\[\norm{r_j}{\shml}:=\sup_{0\neq\psi \in H^{1}(\Om)}\frac{\dual{r(t)}{\psi}}{\normL{\na\psi}}.\]

The following proposition is the counterpart of Proposition \ref{prop2.3}
for the above mixed approximation.

\begin{proposition}\label{prop2.4}
Suppose that $|u_0|, |v_0|\leq 1$, $\eps_0$ and $C_0$ be the same
as in Lemma \ref{lem2.2}. Let $(u,\vphi)$ and $(v,\phi)$ be the
solutions of \eqref{e2.14}-\eqref{e2.16} and
\eqref{e2.17}-\eqref{e2.19}, respectively. Then, for any
$\varepsilon\in (0,\varepsilon_0]$, there exists a positive
constant $C$, which is independent of $\eps$ and $t$, such that
there holds
\begin{equation}\label{e2.20}
\begin{split}
&\normL{\na \De^{-1} (v(t) - u(t))}^2
+\int_0^t\,\Bigl(\eps^4\normL{\nabla(v(s)-u(s))}^2\\
&\hskip 1.9in+\frac{1}{\eps}\norm{v(s)-u(s)}{L^4}^4\Bigr)e^{(2C_0+8)(t-s)}\, ds\\
&\leq C\Bigl[1+\frac{1}{\hat\xi(t)}\Bigr]\int_0^t \left(\norm{r_1}{\shml}^2
+\frac{1}{\eps^2}\norm{r_2}{\shml}^2\right) e^{(2C_0+8)(t-s)}\, ds \\
&\hskip 1.9in
+ \frac{C}{\hat\xi(t)}\normL{\na \De^{-1} (v_0 - u_0)}^2e^{(2C_0+8)t}
\end{split}
\end{equation}
for all $t\in [0,T^{**})$. Here
\begin{equation}\label{e2.20a}
\begin{split}
&\hat\xi(t):=1- C\eps^{-\frac{5(2+N)}{2}} e^{(2C_0+8)t}
\Big\{\normL{\na \De^{-1} (v_0-u_0)}^2  \\
&\hskip 1in
+\int_0^t\left(\norm{r_1}{\shml}^2+\frac{1}{\eps^2}\norm{r_2}{\shml}^2\right)\,
e^{-(2C_0+8)s}\, ds\Big\},
\end{split}
\end{equation}
and $T^{**}\in [0,T]$ satisfying $\hat\xi(T^{**})> 0$.
\end{proposition}

\begin{proof}
Since the proof is very similar to that of Proposition \ref{prop2.3}, we
only highlight the main differences and omit the overlaps.

Let $w(t):=v(t)-u(t)$ and $\theta(t):=\vphi(t)-\phi(t)$.
Subtracting \eqref{e2.14}-\eqref{e2.16} from their corresponding
equations in \eqref{e2.17}-\eqref{e2.19} we get the following
``error" equations: for $t\in [0,T]$
\begin{align}\label{e2.20b}
\bigl(w_t,\psi\bigr)+ \bigl(\na \theta,\na\psi\bigr)
&=\dual{r_1}{\psi} \quad\forall\psi\in H^1(\Om)\, ,\\
\eps \bigl(\na w,\na\chi \bigr)+\frac{1}{\eps} \bigl(f(v)-f(u),
\chi\bigr) -\bigl(\theta,\chi\bigr) &=\dual{\eps r_2}{\chi}
\quad\forall \chi\in H^1(\Om)\, , \label{e2.20c}\\
w(0) &= v_0-u_0\quad\mbox{in } \Om\,. \label{e2.20d}
\end{align}
Setting $\psi=-\De^{-1} w$ in \eqref{e2.20b} and $\chi=w$ in
\eqref{e2.20c} and adding the resulting equations give
\begin{equation}\label{e2.20e}
\begin{split}
&\frac12 \frac{d}{d t} \normL{\na\De^{-1} w}^2 + \frac{1}{\eps}
\norm{w}{L^4}^4 + \eps \normL{\nabla w}^2 +\frac{1}{\eps} \int_\Om
f^\prime(u)\, w^2 dx \\
&\hskip 0.2in
= -\frac{3}{\eps} \int_\Om u\, w^3 dx - \dual{r_1}{\De^{-1} w}
+\dual{\eps r_2}{w}\\
&\hskip 0.2in
\leq \frac{C}{\eps} \int_\Om |w|^3 + \norm{r_1}{\shml}^2+\normL{\na\De^{-1}
w}^2+\frac{1}{\eps^2}\norm{r_2}{\shml}^2+\frac{\eps^4}{4} \normL{\na w}^2.
\end{split}
\end{equation}
Here we have used the identity $\bigl(\theta, w\bigr)
+\bigl(\na\theta,\na \De^{-1}w \bigr)=0$.

Clearly, the only difference between \eqref{e2.20e} and
\eqref{e2.9} is the last four terms  on the right hand side of
\eqref{e2.20e}. Repeating the remaining proof of Proposition
\ref{prop2.3} after \eqref{e2.9}, we see that the conclusion of
Proposition \ref{prop2.3} holds with
$\norm{r_1}{\shml}^2+\frac{1}{\eps^2}\norm{r_2}{\shml}^2$ in the
place of $\eps^{-2}\norm{r}{\shmt}^2$, hence, \eqref{e2.20} holds.
The proof is complete.
\end{proof}

A similar statement to that of Corollary \ref{cor2.1} also holds.
We omit its proof since it is simple.

\begin{corollary}\label{cor2.2}
Under the assumptions of Proposition \ref{prop2.4},
\eqref{e2.20} holds for $T^{**}=T$ if $v_0$ and $(r_1,r_2)$
satisfy the following constraint
\begin{equation*}\label{e2.20f}
\begin{split}
&\left\{\normL{\na \De^{-1} (v_0-u_0)}^2 +\int_0^T
\left(\norm{r_1}{\shml}^2+\frac{1}{\eps^2}\norm{r_2}{\shml}^2\right)
e^{-(2C_0+8)s}\, ds\right\}^{\frac12} \\
&\hskip 1in
\leq C^{-1} e^{-(C_0+4)T} \eps^{\frac{5(2+N)}{4}}
=\left\{
\begin{array}{ll}
        O(\eps^{5}) &\quad\mbox{if } N=2,\\
        O(\eps^{6.25}) &\quad\mbox{if } N=3.
        \end{array} \right.
\end{split}
\end{equation*}
\end{corollary}

We note that Proposition \ref{prop2.4} and Corollary \ref{cor2.2}
only give polynomial order (in $\frac{1}{\eps}$) continuous dependence
estimates for $v-u$. In the next proposition, we derive some
estimates for $\vphi-\phi$.

\begin{proposition}\label{prop2.5}
Under the assumptions of Corollary \ref{cor2.2} there holds
\begin{equation}\label{e2.20j}
\begin{split}
&\int_0^T \norm{\vphi(s)-\phi(s)}{H^{-1}}^{\frac65} ds\\
&\qquad\leq
\frac{C}{\eps^2}\Bigl[1+\frac{1}{\hat\xi(t)}\Bigr]
\int_0^t \left(\norm{r_1}{\shml}^2
+\frac{1}{\eps^2}\norm{r_2}{\shml}^2\right) e^{(2C_0+8)(t-s)}\, ds \\
&\hskip 2.2in
+ \frac{C}{\eps^2\hat\xi(t)}\normL{\na \De^{-1} (v_0 - u_0)}^2 e^{(2C_0+8)t}.
\end{split}
\end{equation}
Moreover, for $N=2$, if $r_2(t)\in L^2(\Omega)$, there also holds
\begin{equation}\label{e2.20i}
\begin{split}
\norm{v(t)-u(t)}{L^2}^2 &+\frac{1}{\eps}\int_0^t
\norm{\vphi(s)-\phi(s)}{L^2}^2 ds\\
 \leq& \frac{C}{\eps^7}\Bigl[1+\frac{1}{\hat\xi(t)}\Bigr] \int_0^t
\left(\norm{r_1}{\shml}^2
+\frac{1}{\eps^2}\norm{r_2}{\shml}^2\right) e^{(2C_0+8)(t-s)}\, ds \\
& + \frac{C}{\eps^7\hat\xi(t)}\normL{\na \De^{-1}(v_0 - u_0)}^2
e^{(2C_0+8)t} + \eps\int_0^t \norm{r_2}{L^2}^2 ds.
\end{split}
\end{equation}
for all $t\in [0,T]$. Where $\hat\xi(t)$ is defined by \eqref{e2.20a}.
\end{proposition}

\begin{proof}
From \eqref{e2.20c}, \eqref{efvu}, and the fact that
$\norm{u}{L^\infty}\leq C$ (cf. \cite{caffarelli95,XA3})  we have for any
$\chi\in H^1_0(\Om)$
\begin{align*}
\bigl(\theta,\chi\bigr) &=\eps \bigl(\na w,\na\chi \bigr)
+\frac{1}{\eps} \bigl(f(v)-f(u), \chi\bigr) -\dual{\eps r_2}{\chi} \\
&\leq \eps \norm{\na w}{L^2} \norm{\na \chi}{L^2}
+\frac{C}{\eps} \Bigl[\, \norm{w}{L^2} \norm{\chi}{L^2}
+\norm{w}{L^{\frac{18}{5}}}^3 \norm{\chi}{L^6}
+ \norm{w}{L^4}^2 \norm{\chi}{L^2}\, \Bigr] \\
&\qquad\qquad + \eps \norm{r_2}{\shml}\norm{\na\chi}{L^2},
\end{align*}
which and the interpolation inequality
\[
\norm{w}{L^{\frac{18}{5}}}\leq
\norm{w}{L^4}^{\frac89} \norm{w}{L^2}^{\frac19}
\]
yield
\begin{align}\label{e2.20k}
\norm{\theta(t)}{H^{-1}}^{\frac65}
&\leq \eps^{\frac65} \norm{\na w(t)}{L^2}^{\frac65}
+C\eps^{-\frac65} \Bigl[\norm{w(t)}{L^2}^2
+ \norm{w(t)}{L^4}^4 \Bigr]
+ \eps^{\frac65}\norm{r_2(t)}{\shml}^{\frac65}.
\end{align}
\eqref{e2.20i} now follows from integrating \eqref{e2.20k} in $t$
over $[0,T]$, and appealing to \eqref{e2.20} and Corollary \ref{cor2.2}.

To show \eqref{e2.20i}, adding \eqref{e2.20b} and \eqref{e2.20c}
after setting $\psi=w$ and $\chi=-\frac{1}{\eps} \theta$, and
using the Schwarz inequality we get
\begin{eqnarray}\nonumber
&&\frac12 \frac{d}{d t}\norm{w}{L^2}^2 +\frac{1}{\eps}\norm{\theta}{L^2}^2
=\dual{r_1}{w} -\dual{r_2}{\theta}
+\frac{1}{\eps^2} \bigl( f(v)-f(u), \theta\bigr) \\
&&\hskip 0.5in
\leq \norm{r_1}{\shml} \norm{\nabla w}{L^2}
+\norm{r_2}{L^2} \norm{\theta}{L^2}
+\frac{C}{\eps^2} \norm{\theta}{L^2} \bigl[\, \norm{w}{L^2}
+\norm{w}{L^6}^3+ \norm{w}{L^4}^2 \,\bigr]
\label{e2.20g} \\
&&\hskip 0.5in
\leq \frac{1}{2\eps} \norm{\theta}{L^2}^2
+ \frac{C}{\eps^3}\Bigr[ \norm{\nabla w}{L^2}^2 %+ \norm{ w}{L^6}^6
+ \norm{ w}{L^4}^4 \Bigr]
+ \norm{r_1}{\shml}^2 +\eps \norm{r_2}{L^2}^2. \nonumber
\end{eqnarray}
Integrating \eqref{e2.20g} over $[0,T]$, the desired estimate
\eqref{e2.20i} follows from an application of \eqref{e2.20}
and Corollary \ref{cor2.2}. The proof is complete.
\end{proof}

\subsection{An abstract framework for a posteriori estimates}\label{sec-2.3}

In this section, we first recall an abstract framework given in
\cite{feng_wu03} for deriving a posteriori estimates based on continuous
dependence estimates of an underlying evolution equation. We refer readers
to a recent survey paper by Cockburn \cite{cockburn03} and the references
therein for applications of a similar method to problems of
hyperbolic conservation law.  We then extend this abstract framework to
mixed approximations of general evolution equations. Since the idea for
deriving a posteriori error estimates essentially works for a large
class of evolution problems, we shall present it in an abstract fashion.

Let $V$ be an Hilbert space and ${\mathcal L}$ be an operator from
$D({\mathcal L})$ ($\subset V$), the domain of ${\mathcal L}$, to $V^*$,
the dual space of $V$. We consider the abstract evolution problem
\begin{alignat}{2}\label{e2.21}
\frac{\partial u}{\partial t} + {\mathcal L}(u)&=r &&\quad\mbox{in }\Om_T, \\
u(0) &=u_0 &&\quad\mbox{in }\Om. \label{e2.22}
\end{alignat}
Suppose that $u^{(j)}$ is the (unique) solution of \eqref{e2.21}-\eqref{e2.22}
with respect to the data $(r^{(j)}, u_0^{(j)})$ for $j=1,2$, respectively.
Assume that $u^{(j)}$ satisfy the continuous dependence estimate
\begin{equation}\label{e2.23}
|||u^{(1)}-u^{(2)}|||\leq F(r^{(1)}-r^{(2)}) + G(u^{(1)}_0-u^{(2)}_0)
\end{equation}
for some (monotone increasing) functionals $F(\cdot)$ and $G(\cdot)$.
Where $|||\cdot|||$ stands for the standard norm in $L^\ell((0,T);V)$ for
some $1\leq \ell\leq \infty$.

The following theorem was proved in \cite{feng_wu03}.

\medskip
\begin{theorem}\label{thm2.1}
Let $u$ denote the solution of \reff{e2.21}-\reff{e2.22}, and $u^A$
be an approximation of $u$ with the initial value $u_0^A$. Suppose
that problem \eqref{e2.21}-\reff{e2.22} satisfies the
continuous dependence estimate \eqref{e2.23}, then there holds
\begin{eqnarray}\label{e2.24}
|||u-u^A||| &\leq& F(R(u^A)) + G(u_0-u^A_0),\\
R(u^A)&:=& r-\frac{\partial u^A}{\partial t} - {\mathcal L}(u^A) .
\label{e2.25}
\end{eqnarray}
\end{theorem}

\begin{remark}\label{rem2.2}
(a). Clearly, the quantity $R(u^A)$ is the residual of $u^A$. This
residual is often difficult to compute or too expensive to compute
exactly. In practice, an upper bound for $R(u^A)$, which should be easy and
cheap to compute, is sought and used to replace $R(u^A)$
in $F(R(u^A))$ in the above a posteriori error estimate. In the next
section we shall give such an estimate for conforming finite element
approximations of the Cahn-Hilliard equation (cf. \cite{ciarlet78,elliottc89}).

(b). A posteriori error estimate \eqref{e2.24} holds for any approximation
$u^A$ of $u$, including non-computable abstract approximations
(cf. \cite{alikakos94}). However, only computable approximations such as
those obtained by finite element methods, finite difference methods,
finite volume methods and spectral methods are of practical interests.
\end{remark}

The above a posteriori estimate can be easily extended to mixed approximations
of problem \eqref{e2.21}-\eqref{e2.22}. We recall that a mixed formulation
of \eqref{e2.21}-\eqref{e2.22} seeks a pair of functions
$(u,p)\in V_1\times V_2$ such that
\begin{eqnarray}\label{e2.26}
\frac{\p u}{\p t} + {\mathcal L}_1 (p) = \mu &&\qquad\mbox{in }\Om_T,\\
p-{\mathcal L}_2 (u) = \eta &&\qquad\mbox{in }\Om_T, \label{e2.27}\\
u(0) = u_0 &&\qquad\mbox{in }\Om. \label{e2.28}
\end{eqnarray}
Where $\{V_i\}_{i=1}^2$ are two Hilbert spaces. ${\mathcal L}_i$
is some operator from $D({\mathcal L}_i)\, (\subset V_i)$, the domain
of ${\mathcal L}_i$, to $V_i^*$, the dual space of $V_i$, which satisfies
${\mathcal L}= {\mathcal L}_1\circ {\mathcal L}_2$. $\mu$ and $\eta$
are two known functions which are appropriately chosen so that
problem \eqref{e2.26}-\eqref{e2.28} is equivalent to problem
\eqref{e2.21}-\eqref{e2.22}.

Suppose that $(u^{(j)},p^{(j)})$ is the (unique) solution of
\eqref{e2.26}-\eqref{e2.28} with respect to the data
$(\mu^{(j)},\eta^{(j)},$ $ u_0^{(j)})$ for $j=1,2$, respectively.
Assume that $(u^{(j)},p^{(j)})$ satisfy the following continuous
dependence estimate
\begin{equation}\label{e2.29}
|||u^{(1)}-u^{(2)}|||_{1} + |||p^{(1)}-p^{(2)}|||_{2}
\leq \Phi(\mu^{(1)}-\mu^{(2)}) + \Psi(\eta^{(1)}-\eta^{(2)}) +
Z(u^{(1)}_0-u^{(2)}_0)
\end{equation}
for some (monotone increasing) nonnegative functionals $\Phi(\cdot)$,
$\Psi(\cdot)$, and $Z(\cdot)$. Where $|||\cdot|||_i$ denotes the
standard norm in $L^\ell((0,T);V_i)$ for some $1\leq \ell\leq \infty$.
Then we have

\medskip
\begin{theorem}\label{thm2.2}
Let $(u,p)$ be the solution of \reff{e2.26}-\reff{e2.28},
and $(u^A,p^A)$ be an approximation of $(u,p)$ with the initial
value $u_0^A$. Suppose that problem \eqref{e2.26}-\reff{e2.28}
satisfies the continuous dependence estimate \eqref{e2.29}, then there
holds
\begin{align}\label{e2.30}
|||u-u^A|||_1 + |||p-p^A|||_2 \leq \Phi(R_1(u^A, p^A)) +
\Psi(R_2(u^A,p^A))
+ Z(u_0-u_0^A) ,&\\
R_1(u^A,p^A):= \mu-\frac{\partial u^A}{\partial t}-{\mathcal
L}_1(p^A) ,\quad R_2(u^A,p^A):= \eta-p^A + {\mathcal L}_2(u^A).&
\label{e2.31}
\end{align}
\end{theorem}

\begin{proof}
Define
\[
\mu^A:= \frac{\partial u^A}{\partial t}+{\mathcal L}_1(p^A),\qquad
\eta^A:= p^A - {\mathcal L}_2(u^A) .
\]
\eqref{e2.30} follows easily from \eqref{e2.29} with $\mu^{(1)}=\mu$,
$\mu^{(2)}=\mu^A$, $\eta^{(1)}=\eta$, $\eta^{(2)}=\eta^A$,
$u_0^{(1)}=u_0$, and $u_0^{(2)}=u_0^A$.
\end{proof}

We conclude this section by the following remark.

\begin{remark}\label{rem2.3}
The quantity $\{R_i(u^A,p^A)\}_{i=1}^2$ are the residuals
of $(u^A,p^A)$, which are often difficult to compute or too expensive to compute
exactly. In practice, an upper bound for $R_i(u^A,P^A)$, which should be easy
and cheap to compute, is sought and used to replace $R_i(u^A,p^A)$
in the terms $\Phi(R_1(u^A,p^A))$ and $\Psi(R_2(u^A,p^A))$ of \eqref{e2.30}.
In the next section we shall give such an estimate for mixed finite element
approximations of the Cahn-Hilliard equation (cf. \cite{elliott89,XA3}).
\end{remark}

\section{A posteriori error estimates for finite element approximations}
\label{sec-3}

In this section we shall apply the abstract frameworks of the
previous section to derive some practical a posteriori error
estimates for conforming finite element approximations of the
Cahn-Hilliard equation and for the Ciarlet-Raviart mixed finite element
approximations of the Cahn-Hilliard equation
\cite{ciarlet78,XA3,scholz78}. As expected, the polynomial order
(in $\frac{1}{\eps}$) continuous dependence estimate of
Propositions \ref{prop2.3} -- \ref{prop2.4} play a critical role.

For $N=2,3$, let $\cT$ be a regular ``triangulation" of $\Om$ such
that $\overline{\Om} = \bigcup_{K \in \cT} \overline{K}$, ($K\in \cT$ are
tetrahedrons in the case $N=3$). Recall that any element $K\in
\cT$ is assumed to be closed.  Let $\cF$ be the set of all faces
(sides in case of $N=2$). For any $K\in\cT$ and $\tau\in\cF$, let
$h_K$ and $h_\tau$ denote the diameters of $K$ and $\tau$,
respectively.

\subsection{Conforming finite element methods}
\label{sec-3.1}

Let $S_h\subset H^2_E(\Om)$ be a conforming finite element
space which consists of  piecewise polynomials on $\cT$ satisfying
the homogeneous Neumann condition. The continuous in time
semi-discrete finite element discretization of
(\ref{e1.1})-(\ref{e1.3}) is defined by seeking $u_h:
[0,T]\rightarrow S_h$ such that for $t\in [0,T]$
\begin{align}\label{euh}
\langle \frac{\p u_h}{\p t}, \psi_h\rangle +\eps\bigl(\De u_h, \De
\psi_h\bigr) +\frac{1}{\eps} \bigl(\na(f(u_h)),\na \psi_h\bigr)
&=0 \quad\forall \psi_h\in  S_h ,
\end{align}
with some starting value $u_h(0)=u_{0h}\in S_h$ satisfying
$\int_\Om u_{0h} dx =\int_\Om u_0 dx$.

For $t\in (0,T]$, we define the residual $r_h(t)\in \bigl(H^2(\Om)
\bigr)^*$ of $u_h$ by
\begin{equation}\label{e3.2}
\langle \frac{\p u_h}{\p t}, \psi\rangle +\eps\bigl(\De u_h, \De
\psi\bigr) +\frac{1}{\eps} \bigl(\na(f(u_h)),\na \psi\bigr)
  = \dual{r_h(t)}{\psi} \quad
\forall \psi \in H^2_E(\Om) .
\end{equation}
Then
\begin{equation}\label{e3.3} \dual{r_h(t)}{\psi_h} = 0\qquad
\forall \psi_h \in S_h .
\end{equation}

\begin{remark}
One can derive a priori error
estimates of $u_h$ which only depends on $\frac{1}{\eps}$ in low
polynomial orders by using the nonstandard analysis of \cite{XA3}.
We refer interested readers to \cite{XA3} for a detailed
exposition.
\end{remark}

It is easy to see that Proposition~\ref{prop2.3},
Proposition~\ref{prop2.2} and Theorem~\ref{thm2.1} all are valid
if both $v$ and $u^A$ are replaced by $u_h$, and both $r$ and $R(u^A)$ are
replaced by $r_h$. Hence, we immediately obtain two a posteriori
error estimates for $u_h-u$. As pointed out in Remark \ref{rem2.2}
(a), for practical considerations, it is necessary to derive an
upper bound for $\norm{r_h}{\shmt}$ which is easy to compute. In
this section we shall establish such a bound, which then leads to
practical a posteriori error estimates for $u_h-u$. To the end, we
need the following local approximation properties of
conforming finite element spaces.

\begin{assumption}\label{lapp}
There exists a interpolant $\Pi_h$ form $H^2_E(\Om)$ to $S_h$ such
that for any $\psi\in H^2_E(\Om)$, $K\in\cT$, and $\tau\in\cF$
\[\norm{\psi-\Pi_h\psi}{L^2(K)}\le C
    h_K^2\norm{\psi}{H^2(\tK)},\]
\[\norm{\psi-\Pi_h\psi}{L^2(\tau)}\le C
    h_\tau^{3/2}\norm{\psi}{H^2(\tilde\tau)},\quad
\norm{\frac{\p(\psi-\Pi_h\psi)}{\p n}}{L^2(\tau)}\le C
    h_\tau^{1/2}\norm{\psi}{H^2(\tilde\tau)},\]
where $C$ is a constant only depending on the minimum angle of the
mesh $\cT$, $\tK$ and $\tilde\tau$ are the union of all elements
having non-empty intersection with $K$ and $\tau$, respectively.
\end{assumption}

\begin{remark}\label{rapp}
It  is not hard to show that Assumption \ref{lapp} is fulfilled by
the well-known confirming elements, including Argyris element and
Bell's element (cf. \cite{ciarlet78}), and the interpolant $\Pi_h$ can
be constructed by following the idea of Scott-Zhang interpolation \cite{sz90}.
\end{remark}

For any $K\in\cT$, introduce the element residual
\begin{equation}\label{eres}
    R_K(t) = \frac{\p u_h(t)|_K}{\p t}+
    \De \big(\eps\De u_h(t)|_K-\frac1{\eps}f(u_h(t)|_K)\big).
\end{equation}
For any face $\tau\in\cF$ of element $K$ we define two kinds of
residual jumps across $\tau$. If $\tau$ is an interior face which
is the common face between $K$ and $K^\prime$, let
\begin{equation}\label{ejumpi}
    J_\tau(t)=\big(\na\De u_h(t)|_{K^\prime}-\na\De u_h(t)|_{K}\big)\cdot
    n,\quad \Hat{J}_\tau(t)=\De u_h(t)|_{K}-\De u_h(t)|_{K^\prime}\,.
\end{equation}
Here $n$ denotes the unit outer normal vector to $\tau$. If
$\tau\subset\p\Om$ is a boundary face, define
\begin{equation}\label{ejumpb}
    J_\tau(t)=-2\na\De u_h(t)|_{K}\cdot
    n,\qquad \Hat{J}_\tau(t)=2\De u_h(t)|_{K}\,.
\end{equation}
For any $K\in\cT$, let $\eta_K$ denote the following local error estimator
\begin{equation}\label{eeta}
    \eta_K(t)=h_K^2\norm{R_K}{L^2(K)}+\sum_{\tau\subset\p
    K}\left(\frac{h_\tau^3}{2}\norm{J_\tau}{L^2(\tau)}^2 +
    \frac{h_\tau}{2}\norm{\Hat{J}_\tau}{L^2(\tau)}^2\right)^{1/2} .
\end{equation}
Next we estimate the residual $r_h(t)$ in terms of $\eta_K(t)$.

\begin{proposition}\label{prop3.1}
There exists a constant $C$, which depends only on the minimum angle
of the mesh $\cT$, such that
\begin{equation}\label{eest_r}
    \norm{r_h(t)}{\shmlt}^2\le
    C\sum_{K\in\cT}\big(\eta_K(t)\big)^2 .
\end{equation}
\end{proposition}

\begin{proof}
By \eqref{e3.2}, \eqref{e3.3}, and integration by parts we obtain
for any $\psi\in H^2_E(\Om)$ and $\psi_h\in S_h$
\begin{align*}
&\dual{r_h(t)}{\psi}=\dual{r_h(t)}{\psi-\psi_h}\\
&\quad =\langle \frac{\p u_h}{\p t}, \psi-\psi_h\rangle +\eps\bigl(\De u_h,
\De(\psi-\psi_h) \bigr) +\frac{1}{\eps} \bigl(\na(f(u_h)),\na
(\psi-\psi_h)\bigr)\\
&\quad =\sum_{K\in\cT}\Big\{\int_K\Big(\frac{\p u_h}{\p
t}+\De\big(\eps\De u_h-\frac1{\eps}f(u_h)\big)\Big)(\psi-\psi_h) dx\\
&\quad\, +\int_{\p K}\Big(-\frac{\p\De u_h}{\p n}(\psi-\psi_h)+\De
u_h\frac{\p(\psi-\psi_h)}{\p n}\Big) d\sigma
+ \int_{\p K}\frac1\eps\frac{\p f(u_h)}{\p
n}(\psi-\psi_h) d\sigma \Big\} .
\end{align*}
Since any interior face be a common face of two elements whose
outer normal vectors to the face are opposite in direction, on noting
that $u_h\in C^1$ we get
\begin{align*}
   \dual{r_h(t)}{\psi} =\sum_{K\in\cT}\Big\{&\int_K R_K(\psi-\psi_h) dx\\
&+\frac12\sum_{\tau\subset \p K}\int_{\p
K}\Big(J_\tau(t)(\psi-\psi_h)+\Hat{J}_\tau(t)\frac{\p(\psi-\psi_h)}{\p
n}\Big) d\sigma \Big\}\,.
\end{align*}
Choosing $\psi_h=\Pi_h\psi$, the desired estimate \eqref{eest_r}
follows from an application of the Schwarz inequality and
Assumption \ref{lapp}. The proof is complete.
\end{proof}

Combining Proposition \ref{prop3.1}, \ref{prop2.3}--\ref{prop2.2},
and Corollary \ref{cor2.1}, we immediately obtain the following
theorem which presents a posteriori error estimates for the
finite element method.

\begin{theorem}\label{t3.1}
Suppose that $|u_0|, |u_{0h}|\leq 1$, and that
$\int_\Om(u_0-u_{0h}) dx = 0$. Let $\eps_0$ and $C_0$ be the same as
in Lemma \ref{lem2.2}, $u$ and $u_h$ be the solutions of
\eqref{eu}-\eqref{eu1} and \eqref{euh}, respectively. Define
\begin{equation}\label{e3.4}
\begin{split}
\xi_h(t):=1- C&\eps^{-\frac{5(2+N)}{2}} \, e^{(2C_0+8)t}
\times\\
 &\left\{\normL{\na \De^{-1} (u_{0h}-u_0)}^2
+\frac{1}{\eps^2} \int_0^t\, e^{-(2C_0+8)s}\,
\sum_{K\in\cT}\eta_K^2(s)\, ds\right\}\,.
\end{split}
\end{equation}
Assume $\xi_h(T)>0$. Then, for any $\varepsilon\in
(0,\varepsilon_0]$ and $t\in [0,T]$, the following a posteriori
error estimates hold
\begin{align}\label{e3.5}
\begin{split} &\normL{\na \De^{-1} (u_h(t) - u(t))}^2
+\int_0^t\,\Bigl(\eps^4\normL{\nabla(u_h(s)-u(s))}^2 \\ 
&\hskip 1.9in
+ \frac{1}{\eps}\norm{u_h(s)-u(s)}{L^4}^4\Bigr)e^{(2C_0+8)(t-s)}\, ds\\
&
\leq \xi_h(t)^{-1}\normL{\na \De^{-1} (u_{0h} - u_0)}^2e^{(2C_0+8)t} \\
&\hskip 1in
+\left[1+\frac{1}{\xi_h(t)}\right]\frac{C}{\eps^2}
\int_0^t e^{(2C_0+8)(t-s)}\, \sum_{K\in\cT}\eta_K^2(s)\, ds .
\end{split}
\end{align}
\begin{align}\label{e3.6}
\begin{split} &\normL{u_h(t) - u(t)}^2 +\int_0^t\Big(\eps\normL{\De
(u_h(s)-u(s))}^2\\
&\hskip 1.6in +\frac{6}{\eps}\normL{(u_h(s)-u(s))\na
(u_h(s)-u(s))}^2\Big)\,ds
\\
&\leq \normL{u_{0h}-u_0}^2+\frac{C}{\eps^5\xi_h(t)}\normL{\na\De^{-1}(u_{0h}
- u_0)}^2e^{(2C_0+8)t}  \\
&\hskip 0.96in
+\frac{C}{\eps^7}\left[1+\frac{1}{\xi_h(t)}\right]\int_0^t\, e^{(2C_0+8)(t-s)}
\sum_{K\in\cT}\eta_K^2(s)\, ds\,.
\end{split}
\end{align}
\end{theorem}

\subsection{Ciarlet-Raviart mixed finite element methods}\label{sec-3.1.2}

Let $V_h^m$ denote the $P_m\, (m\geq 1)$ conforming finite element
subspace of $H^1(\Om)$ consisting of {\em continuous} piecewise
$m^{\textup{th}}$ order polynomial functions on $\cT$ (cf.
\cite{ciarlet78}), that is,
\begin{equation}\label{e3.20}
V_h^m = \Bigl\{v_h \in C(\overline{\Om});\, v_h\bigl \vert_K\in
P_m(K)\quad \forall\, K \in \cT \Big\}.
\end{equation}

Following \cite{elliott89, XA3}, the continuous in time
semi-discrete mixed finite element method is defined to find
$(u_h,\vphi_h): [0,t]\rightarrow [V_h^m]^2$ such that for $t\in
(0,T]$
\begin{align} \label{e3.21}
\Bigl(\frac{\p u_h}{\p t}, \psi_h\Bigr) + \bigl(\nabla \vphi_h,
\nabla \psi_h) &=0 \qquad
\forall\, \psi_h \in V_h^m\, , \\
\eps\, \bigl( \nabla u_h, \nabla \chi_h\bigr) + \frac{1}{\eps}\,
\bigl(f(u_h), \chi_h\bigr) -\bigl(\vphi_h, \chi_h\bigr) &= 0
\qquad \forall\, \chi_h \in V_h^m\, ,  \label{e3.22}
\end{align}
with some suitable starting value $u_h(0)=u_{0h}\in V_h^m$
satisfying $\int_\Om u_{0h} dx =\int_\Om u_0 dx $.

We remark that the finite element spaces $V_h^m\times V_h^m$
is a family of stable mixed finite spaces known as the
Ciarlet-Raviart mixed finite elements for the biharmonic problem
(cf. \cite{ciarlet78,scholz78}), that means the following {\em
inf-sup} condition holds
\begin{equation}\label{e3.23}
\inf_{0\not\equiv\chi_h\in V_h^m}\, \sup_{0\not\equiv\psi_h\in
V_h^m} \frac{(\nabla \psi_h, \nabla
\chi_h)}{\norm{\psi_h}{H^1}\,\norm{\chi_h}{H^1}} \geq c_0
\end{equation}
for some $h$-independent constant $c_0 >0$.

We also define the residual $(\mu_h(t),\eta_h(t))\in [\shml]^2$ of
$(u_h,\vphi_h)$ by
\begin{align}\label{e3.24}
\Bigl(\frac{\p u_h}{\p t}, \psi\Bigr) + \bigl(\nabla \vphi_h,
\nabla \psi) =\dual{r_h^{(1)}(t)}{\psi} &\qquad
\forall\, \psi \in H^1(\Om)\, , \\
\eps\, \bigl( \nabla u_h, \nabla \chi\bigr) + \frac{1}{\eps}\,
\bigl(f(u_h), \chi\bigr) -\bigl(\vphi_h, \chi\bigr) = \dual{\eps
r_h^{(2)}(t)}{\chi} &\qquad \forall\, \chi \in H^1(\Om), .
\label{e3.25}
\end{align}
Clearly, there holds
\begin{equation}\label{e3.26}
\dual{r_h^{(1)}(t)}{\psi_h}=\dual{r_h^{(2)}(t)}{\chi_h}=0 \qquad\forall
(\psi_h,\chi_h)\in [V^m_h]^2 .
\end{equation}

For any $K\in \cT$, we introduce the element residual
\begin{equation}\label{erese}
\begin{split}
R_K^{(1)}(t) &:= \frac{d u_h(t)|_K}{d t} - \De
\bigl(\vphi_h(t)|_K\bigr) \, ,\\
R_K^{(2)}(t) &:= -\De \bigl(u_h(t)|_K\bigr)+\frac1{\eps^2}
f(u_h(t)|_K)-\frac1\eps\vphi_h(t)\,.
\end{split}
\end{equation}
For any common face $\tau$ of $K_1$
$K_2\in \cT$, we define the residual jumps across $\tau$ as
\begin{equation}\label{eresj}
\begin{split}
J_\tau^{(1)}(t)&=\left(\na
\vphi_h(t)|_{K_1}-\na\vphi_h(t)|_{K_2}\right)\cdot n_1,\\
J_\tau^{(2)}(t)&=\left(\na u_h(t)|_{K_1}-\na u_h(t)|_{K_2}\right)\cdot n_1 ,
\end{split}
\end{equation}
where $n_1$ is the unit normal vector to $\tau$ pointing from $K_1$
to $K_2$. For any $\tau\subset\p\Om$ which is a face of some element $K$,
let\begin{equation}\label{eresj1}
J_\tau^{(1)}(t)=2\na\vphi_h(t)|_{K}\cdot n,\quad
 J_\tau^{(2)}(t)=2\na u_h(t)|_{K}\cdot n\, .
\end{equation}
 For any $K\in\cT$, define the local error estimators with
respect to $K$ as follows
\begin{equation}\label{eeta12}
    \eta_K^{(j)}(t)=h_K\norm{R_K^{(j)}}{L^2(K)}+\sum_{\tau\subset\p
    K}\left(\frac{1}{2}h_\tau\norm{J_\tau^{(j)}}{L^2(\tau)}^2\right)^{\frac12} ,
    \quad j=1,2.
\end{equation}

\begin{proposition}\label{prop3.2}
The following estimate holds for the residual $r_h^{(j)}(t)$
\begin{align}\label{epost}
\norm{r_h^{(j)}(t)}{\shml}^2 &\leq
C\sum_{K\in\cT}\big(\eta_{K}^{(j)}(t)\big)^2,\qquad j=1,2,
\end{align}
where $C$ is some constant which depends only on the minimum
angle of the mesh $\cT$.
\end{proposition}

\begin{proof} By \eqref{e3.24}--\eqref{e3.26} and integration by
parts we obtain that for any $\psi, \chi\in H^1(\Om)$ and $\psi_h,
\chi_h\in V_h^m$
\begin{align*}
\dual{r_h^{(1)}(t)}{\psi}&=\dual{r_h^{(1)}(t)}{\psi-\psi_h}=\Bigl(\frac{\p
u_h}{\p t}, \psi-\psi_h\Bigr) + \bigl(\nabla \vphi_h, \nabla
(\psi-\psi_h)\big)\\
 &=\sum_{K\in\cT}\left(\int_K(u_{h t}-\De\vphi_h)(\psi-\psi_h) dx
      +\int_{\p K}\frac{\p\vphi_h}{\p n}(\psi-\psi_h) d\sigma \right) .
\\
\dual{\eps\, r_h^{(2)}(t)}{\chi} &=\eps\, \bigl( \nabla u_h, \nabla
(\chi-\chi_h)\bigr) + \frac{1}{\eps}\,
\bigl(f(u_h), \chi-\chi_h\bigr) -\bigl(\vphi_h, \chi-\chi_h\bigr)\\
 =\sum_{K\in\cT}&\left(\int_K\big(-\eps\De u_{h}+\frac1\eps f(u_h)-\vphi_h\big)(\chi-\chi_h) dx
      +\eps\int_{\p K}\frac{\p u_h}{\p n}(\chi-\chi_h) d\sigma \right) .
\end{align*}
From the definitions \eqref{erese}--\eqref{eeta12}, we conclude
that
\begin{equation}\label{erj}
\dual{r_h^{(j)}(t)}{\psi}=\sum_{K\in\cT}\left(\int_K
R_K^{(j)}(t)(\psi-\psi_h)
+\frac12\sum_{\tau\subset\p K}\int_{\tau} J_\tau^{(j)}(t)(\psi-\psi_h)\right) .
\end{equation}
Choosing $\psi_h=\Pi_h\psi$, where $\Pi_h$ is the Scott-Zhang
interpolant \cite{sz90}, then the desired estimate \eqref{epost}
follows from an application of the Schwarz inequality and
following approximation properties of the Scott-Zhang
interpolation
\[\norm{\psi-\Pi_h\psi}{L^2(K)}\le C
    h_K\norm{\psi}{H^1(\tK)},\qquad\norm{\psi-\Pi_h\psi}{L^2(\tau)}\le C
    h_\tau^{1/2}\norm{\psi}{H^1(\tilde\tau)}\]
where $C$ is a constant only depending on the minimum angle of the
mesh $\cT$, $\tK$ and $\tilde\tau$ are the union of all elements
having non-empty intersection with $K$ and $\tau$, respectively.
The proof is complete.
\end{proof}

Combining Proposition \ref{prop3.2}, \ref{prop2.4} and Corollary
\ref{cor2.2}, we immediately obtain the following theorem which
presents a posteriori error estimates for the mixed finite element methods.

\begin{theorem}\label{thm3.2}
Suppose that $|u_0|, |u_{0h}|\leq 1$, and that
$\int_\Om(u_0-u_{0h}) dx = 0$. Let $\eps_0$ and $C_0$ be the same as
in Lemma \ref{lem2.2}, and $(u,\vphi)$ and $(u_h,\vphi_h)$ be the
solutions of \eqref{e2.14}-\eqref{e2.16} and
\eqref{e3.21}-\eqref{e3.22}, respectively. Define
\begin{equation}\label{eetamix}
    \eta_K(t)=\left(\big(\eta_K^{(1)}(t)\big)^2
    +\frac1{\eps^2}\big(\eta_K^{(2)}(t)\big)^2\right)^{\frac12},
\end{equation}
and
\begin{equation}\label{e3.28}
\begin{split}
\hat\xi_h(t):=1-C&\eps^{-\frac{5(2+N)}{2}} \,
e^{(2C_0+8)t}\times\\ 
&\Big\{\normL{\na \De^{-1} (v_0-u_0)}^2 +\int_0^t\,
e^{-(2C_0+8)s} \sum_{K\in\cT}\big(\eta_K(s)\big)^2\, ds\Big\}\,.
\end{split}
\end{equation}
Assume $\hat\xi_h(T)> 0$. Then, for any $\varepsilon\in
(0,\varepsilon_0]$, there hold
\begin{align}\label{e3.27}
\begin{split}
&\normL{\na \De^{-1} (u_h(t) - u(t))}^2 +\int_0^t\,\Bigl(\eps^4\normL{\nabla(u_h(s)-u(s))}^2\\
&\hskip 2in+\frac{1}{\eps}\norm{u_h(s)-u(s)}{L^4}^4\Bigr)e^{(2C_0+8)(t-s)}\, ds\\
&
\leq \hat\xi_h(t)^{-1}\normL{\na \De^{-1} (u_{0h} - u_0)}^2 e^{(2C_0+8)t} \\
&\hskip 1.5in
+C\Bigl[ 1+\frac{1}{\hat\xi_h(t)} \Bigr]\int_0^t
e^{(2C_0+8)(t-s)} \sum_{K\in\cT}\big(\eta_K(s)\big)^2\, ds ,
\end{split}
\end{align}
\begin{align}\label{e3.29}
\begin{split} &\int_0^T
\norm{\vphi_h(s)-\vphi(s)}{H^{-1}}^{\frac65} ds \leq
\frac{C}{\eps^2}\Bigl[1+\frac{1}{\hat\xi_h(t)}\Bigr]
\int_0^t e^{(2C_0+8)(t-s)} \sum_{K\in\cT}\big(\eta_K(s)\big)^2 \, ds  \\
&\hskip 2in
+\frac{C}{\eps^2\hat\xi_h(t)}\normL{\na\De^{-1} (u_{0h} - u_0)}^2 e^{(2C_0+8)t}
\end{split}
\end{align}
for all $t\in [0,T)$.
\end{theorem}

\section{An adaptive algorithm}\label{sec-4}

We now present an adaptive algorithm based on the technique of
``method of lines" \cite{bfs01}, i.e., we use the stiff ODE solver
of NDF \cite{sr97} which is a modification of BDF for temporal
integration, and the conforming Argyris element for spatial discretization. The temporal errors are controlled by NDF and assumed
to be sufficiently small that we concentrate solely on controlling
spatial discretization errors. Our local a posteriori error
estimates (cf. Proposition~\ref{prop3.1}) are used to refine and coarsen
the meshes locally. The following adaptive algorithm is an improvement
of the one proposed in \cite{wu02} and is more suitable for computing
the solution of the Cahn-Hilliard equation, which is smooth but contains
a sharp moving front.

\bigskip
\noindent
{\bf Algorithm 4.1.}

For a given tolerance $TOL$, perform the following steps:
\begin{enumerate}
  \item[(i)] Determine an initial mesh $\mathcal{T}_0$ and initial
approximation $u_h(0)$  such that
  $\abs{ u_h(0)- u(0)}_{H^2}$ $ < TOL \times\max(\abs{u_h(0)}_{H^2},1)$. Set  $i=0$.
  \item[(ii)] Do temporal integration $N (=15)$ steps. Denote by $t_{i+1}$ the
current time, and by $n_i$ the number of elements in $K_i$.
  \item[(iii)] Calculate the posteriori error estimate at $t_{i+1}$ :
\[
E_{i+1} = \left(\sum_{j=1}^{n_i}\tilde\eta_{K_j}^2\right)^{1/2},
\quad \tilde\eta_{K_j} = \eta_{K_j}/\max(\abs{u_h(t_{i+1})}_{H^2},1) .
\]
Assume that
$\tilde\eta_{K_1}\leq\tilde\eta_{K_2}\leq\cdots\leq\tilde\eta_{K_{n_i}}$.
\item[(iv)] If  $E_{i+1}>TOL$, then choose  $nr$ such that
\[
nr=\min\left\{j;\,\tilde\eta_{K_j}\geq\frac{1}{2}\tilde\eta_{K_{n_i}},\,
\sum_{l=j}^{n_i}\tilde\eta_{K_l}^2\leq\frac{4}{3}
\left(E_{i+1}^2-TOL^2\right) \right\}.
\]
And refine elements $K_{nr},\cdots,K_{n_i}$ to obtain a new mesh
denoted also by $\mathcal{T}_i$. Redo temporal integration from $t_i$ to
$t_{i+1}$ on the finer mesh. Then go to (iii).
\item[(v)] If  $E_{i+1}\leq TOL$, then choose  $nc$ such that
\[
nc=\max\left\{j;\,\sum_{l=1}^{j}\tilde\eta_{K_l}^2\leq\frac{1}{255}
\left(TOL^2-E_{i+1}^2\right)\right\}.
\]
And coarsen elements $K_1,\cdots,K_{nc}$ to obtain a new mesh
denoted  by $\mathcal{T}_{i+1}.$ Set $i=i+1,$  go to (ii).
\end{enumerate}

In Section~\ref{sec-6}, we shall provide some numerical tests
to gauge performance of the above adaptive algorithm and our a
posteriori error estimates. Our numerical tests show that the algorithm
and the error estimators work remarkably well for the Cahn-Hilliard
equation.

\section{Approximation of the Hele-Shaw flow} \label{sec-5}

Let $\{\Gamma_t^\eps\}_{t\geq 0}$ denote the zero level sets of the
solution $u^\eps$ to the Cahn-Hilliard problem
\eqref{e1.1}-\eqref{e1.3}, and $\{\Gamma_t^{\eps,h}\}_{t\geq 0}$ denote
the zero level sets of the numerical solution $u^\eps_h$ to the
scheme \eqref{euh}. %or to the scheme \eqref{e3.21}-\eqref{e3.22}.
Note that we have put back the super-index $\eps$ on both
$u^\eps$ and $u^\eps_h$ in this section. An interesting (and hard) problem
is to establish the convergence of the numerical
interface $\Gamma_t^{\eps,h}$ to the
true interface $\Gamma_t$ of the Hele-Shaw problem, and also to derive
an a posteriori error estimate for them.  In the following we shall explain
that this can be done in a similar way to that used to derive a
priori error estimates for the numerical interface in \cite{XA3}.

As for all phase field models, the convergence of the numerical
interface to the interface of the limiting problem is usually proved
in two steps.  First, one establishes the convergence of $\Gamma_t^\eps$
to $\Gamma_t$, Second, one proves the convergence of $\Gamma_t^{\eps,h}$
to $\Gamma_t^\eps$.  A triangle inequality then immediately implies the
convergence of $\Gamma_t^{\eps,h}$ to $\Gamma_t$.

For the Cahn-Hilliard equation, we recall that the required
first step was already proved in \cite{alikakos94}. In particular,
we cite the following theorem of \cite{alikakos94}.

\begin{theorem}\label{thm5.0}
Let $\Omega$ be a given smooth domain and $\Gamma_{00}$ be a smooth
closed hypersurface in $\Omega$. Suppose that the Hele-Shaw problem
(\ref{eqn-1.6})-(\ref{eqn-1.10}) starting from $\Gamma_{00}$ has
a smooth solution $(w,\Gamma:=\cup_{0\leq t\leq T}(\Gamma_t\times\{t\}))$ in
the time interval $[0,T]$ such that $\Gamma_t\subset \Omega$ for all
$t\in [0,T]$. Then there exists a family of smooth functions
$\{u^\varepsilon_0(x)\}_{0<\varepsilon \leq 1}$ which are
uniformly bounded in $\eps\in (0,1]$ and $(x,t)\in \overline{\Om}_T$,
such that if $u^\eps$ solves the Cahn-Hilliard equation
(\ref{e1.1})-(\ref{e1.3}) with the initial condition $u^\eps(\cdot,t)
=u^\eps_0(\cdot)$, then
\begin{eqnarray*}
&{\rm (i)}& \quad \lim_{\eps\rightarrow 0} u^\eps(x,t)
=\left\{ \begin{array}{ll}
          1 &\quad\mbox{if } (x,t)\in \mathcal{O}\\
          -1 &\quad\mbox{if } (x,t)\in \mathcal{I}\\
         \end{array} \right.
\, \mbox{uniformly on compact subsets},\\
&{\rm (ii)}& \quad\lim_{\eps\rightarrow 0} \bigl(\frac{1}{\eps} f(u^\eps)
- \eps \Delta u^\eps \bigr) (x,t) = w(x,t) \quad\mbox{uniformly on }
\overline{\Om}_T\, .
\end{eqnarray*}
Where
\[
\mathcal{I}:=\{(x,t)\in \Om\times [0,T]\, ;\, d(x,t)<0\},\qquad
\mathcal{O}:=\{(x,t)\in \Om\times [0,T]\, ;\, d(x,t)>0\}\, ,
\]
and $d(x,t)$ denotes the signed distance function to $\Gamma_t$.
\end{theorem}

Next, we shall prove an a posteriori convergence result for the distance
between $\{\Gamma_t\}_{t\geq 0}$ and $\{\Gamma_t^{\eps,h}\}_{t\geq 0}$,
in particular, the estimate allows one to adjust the mesh size $h$
such that this distance is as small as one wishes before the onset
of singularities.

\begin{theorem}\label{thm5.1}
Let $t_*$ denote the first time when the classical solution of
the Hele-Shaw problem has a singularity. Suppose that
$\Gamma_0=\{x\in\overline{\Omega}; u^\eps_0(x)=0 \}$ is a smooth hypersurface
compactly contained in $\Omega$, and let $\zeta_h(t)$ be same as in
Theorem \ref{t3.1}. Then, for any $\delta\in (0,1)$, there exists
a constant $\hat{\eps}_0 >0$ such that for $t< t_*$
\begin{equation*}\label{e5.4}
\sup_{x\in \Gamma^{\varepsilon,h}_t}
\{\,\mbox{\rm dist}(x, \Gamma_t)\, \} \leq \delta
\quad\qquad\forall \varepsilon\in (0,\hat{\eps}_0) \qquad
\mbox{uniformly on } [0,T],
\end{equation*}
provided that the mesh size $h$ and the starting value $u_h(0)$ satisfy
\begin{eqnarray}\label{e5.3a}
\norm{I_h u^\eps -u^\eps}{L^\infty} &<& \frac{\delta}4,\\
h^{-\frac{N}{2}} \Bigl\{ \norm{u^\eps_0-u^\eps_h(0)}{L^2}
+ \frac{C}{\sqrt{\eps^5\zeta_h(t)}} e^{(4+C_0)T}
\norm{ \na\Delta^{-1}(u^\eps_0-u^\eps_h(0)) }{L^2}\Bigr\}
 &<& \frac{\delta}4, \label{e5.3b} \\
h^{-\frac{N}{2}}  \left\{ \frac{C}{\eps^7} \Bigl[ 1+ \frac{1}{\zeta_h(t)}
\Bigr] \int_0^T  e^{(2C_0+8)(t-s)} \sum_{K\in {\mathcal T}_h}
\eta_K^2(s) ds \right\}^{\frac12} &<& \frac{\delta}4 , \label{e5.3c}
\end{eqnarray}
where $I_h$ denotes standard nodal interpolation operator into
the finite element space $S_h$ (cf. \cite{ciarlet78}).
\end{theorem}

\begin{proof}
First, we prove that $u^\eps_h$ converges uniformly to $1$ on every
compact subset of $\mathcal{O}$.  Let $A$ be a compact subset of
$\mathcal{O}$, for any $(x,t)\in A$, by the triangle inequality we get
\begin{equation}\label{e5.5}
|u^\eps_h(x,t)-1|\leq \norm{u^\eps_h-u^\eps}{L^\infty} + |u^\eps-1|.
\end{equation}

It follows from the inverse inequality, Theorem \ref{t3.1}, and the
assumptions \eqref{e5.3a}--\eqref{e5.3c} that
\begin{align}\label{e5.6}
\norm{u^\eps_h-u^\eps}{L^\infty} &\leq \norm{u^\eps_h-I_h u^\eps}{L^\infty}
+\norm{I_h u^\eps -u^\eps}{L^\infty} \\
&\leq h^{-\frac{N}{2}} \bigl\{ \norm{u^\eps_h- u^\eps}{L^2}
+\norm{u^\eps-I_h u^\eps}{L^2} \bigr\} + \norm{I_h u^\eps -u^\eps}{L^\infty}
\leq \frac{3\delta}4, \nonumber
\end{align}
which together with \eqref{e5.5}, and Theorem \ref{thm5.0}
imply that there exists $\eps_0> 0$ such that
\begin{equation}\label{e5.7}
|u^\eps_h(x,t)-1|\leq \delta\qquad \forall \eps\in (0,\eps_0),\quad (x,t)\in A.
\end{equation}

Similarly, we can show that $u^\eps_h$ converges uniformly to $(-1)$ on every
compact subset of $\mathcal{I}$, that is, there exists $\hat{\eps}_0\in
(0,\eps_0)$ such that for any compact subset $B$ of $\mathcal{I}$ there holds
\begin{equation}\label{e5.8}
|u^\eps_h(x,t)+1|\leq \delta  \qquad \forall \eps\in (0,\hat{\eps}_0),\quad
 (x,t)\in B.
\end{equation}

Define the (open) tabular neighborhood $\mathcal{N}_{\delta}$ of
width $2 \delta$ of $\Gamma_t$ as
\begin{equation}\label{e5.9}
\mathcal{N}_{\delta}:=\{\, (x,t)\in \Om_T\,;\,  d(x,t)< \delta\,\}\, .
\end{equation}
Let $A$ and $B$ now denote the complements of $\mathcal{N}_{\delta}$
in $\mathcal{O}$ and $\mathcal{I}$, respectively, that is,
\[
A=\mathcal{O}\setminus \mathcal{N}_{\delta}\,, \qquad
B=\mathcal{I}\setminus \mathcal{N}_{\delta}.
\]

Note that $A$ is a compact subset of $\mathcal{O}$ and $B$ is a compact subset
of $\mathcal{I}$. Hence, it follows from \eqref{e5.7} and \eqref{e5.8} that
for any $\eps\in (0, \hat{\eps}_0)$
\begin{eqnarray}\label{e5.10}
|u^\eps_h(x,t)-1| &\leq& \delta \qquad\forall\, (x,t)\in A\, ,\\
|u^\eps_h(x,t)+1| &\leq& \delta \qquad\forall\, (x,t)\in B\, .
\label{e5.11}
\end{eqnarray}

Now for any $t\in [0,T]$ and $x\in \Gamma^{\eps,h}_t$, since $u^\eps_h(x,t)=0$, we have
\begin{eqnarray}\label{e5.12}
&&|u^\eps_h(x,t)-1| = 1\, , \\
&&|u^\eps_h(x,t)+1| = 1\, . \label{e5.13}
\end{eqnarray}

Evidently, (\ref{e5.10}) and (\ref{e5.12}) imply that
$(x,t)\not\in A$, and (\ref{e5.11}) and (\ref{e5.13})
says that $(x,t)\not\in B$. Hence $(x,t)$ must reside in the
tubular neighborhood $\mathcal{N}_{\delta}$. Since $t$ is an arbitrary number
in $[0,T]$ and $x$ is an arbitrary point on $\Gamma^{\eps,h}_t$,
therefore, for any $\eps\in (0, \hat{\eps}_0)$
\begin{equation}\label{e5.14}
\sup_{x\in \Gamma^{\eps,h}_t} \bigl( \mbox{\rm dist}(x,\Gamma_t) \bigr)
\leq \delta \quad\mbox{uniformly on } [0,T]\, .
\end{equation}
The proof is complete.

\end{proof}

\section{Numerical Experiments}\label{sec-6}

We shall present a few numerical tests in this section to gauge the
performance of the proposed adaptive algorithm and a posteriori
error estimators. These tests indicate that the algorithm works very
well for the Cahn-Hilliard equation. In all tests to be
given in the following, we take $\Om=[-1,1]^2$.

\bigskip
{\bf Test 1:} Consider the Cahn-Hilliard equation
\eqref{e1.1}-\eqref{e1.3} with the following initial condition
\begin{equation}\label{E1}
\begin{split}
u_0(x,y)= \tanh\bigl(((x-0.3)^2+y^2-0.25^2)/\eps\bigr)
 \tanh\bigl(((x+0.3)^2+y^2-0.3^2)/\eps\bigr).
\end{split}
\end{equation}
Here $\tanh(x)=\dfrac{e^x-e^{-x}}{e^x+e^{-x}}$.

\begin{figure}[hbt]
\centerline{
\includegraphics[scale=0.28]{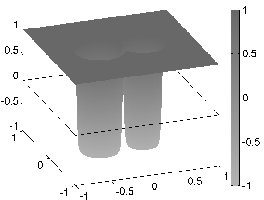}
\includegraphics[scale=0.28]{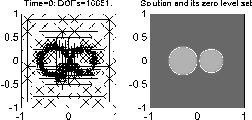}
} \caption{The profile of $u_0$ and its zero level set of Test
1}\label{fig1}
\end{figure}

\begin{figure}[htb]
\centerline{
\includegraphics[scale=0.27]{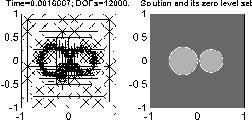}
\includegraphics[scale=0.27]{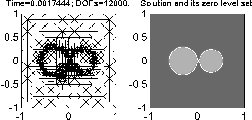}
} \centerline{
\includegraphics[scale=0.27]{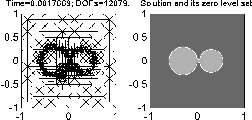}
\includegraphics[scale=0.27]{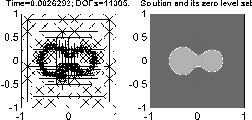}
} \centerline{
\includegraphics[scale=0.27]{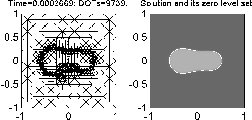}
\includegraphics[scale=0.27]{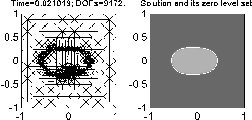}
} \centerline{
\includegraphics[scale=0.27]{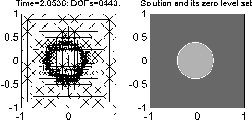}
\includegraphics[scale=0.27]{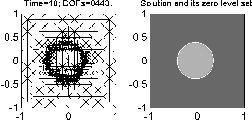}
}  \caption{Snapshots of computed solutions and adaptive meshes for
Test 1}\label{fig2}
\end{figure}

Figure~\ref{fig1} displays the graph of the initial function $u_0$
and its zero level set, which encloses two circles with radii $0.25$
and $0.3$, respectively. It also shows the initial mesh and computed
initial zero level set $\Gamma_0^{0.01,h}$. Figure~\ref{fig2} shows
snapshots of the solution (and its zero level set) of the
Cahn-Hilliard equation and the (adaptive) mesh on which the solution
is computed at $15$ different time steps. $\eps=0.01$ and $TOL=0.02$
are used in the simulation. As expected, the fine mesh follows the
zero level set as it moves. We also note that the number of elements
in the initial mesh $\mathcal{T}_0$ is $3,674$, the minimum area of
the elements is $1.5259\times 10^{-5}$. If a uniform mesh is used,
we need $\frac{4}{1.5259}\times 10^{5}\approx 262,140$ elements and
about $1,180,000$ DOFs.

\begin{figure}[htb]
\centerline{
\includegraphics[scale=0.18]{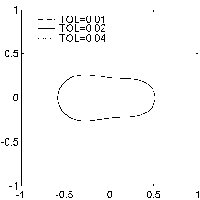}
\includegraphics[scale=0.18]{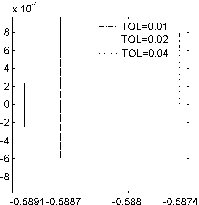}
\includegraphics[scale=0.18]{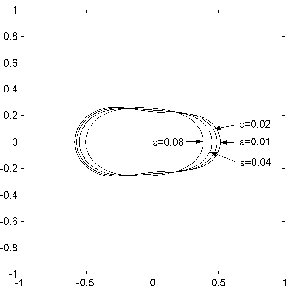}
} \hskip 0.8truein (a) \hskip 1.5truein (b) \hskip 1.5truein (c)
\caption{Convergence of numerical interface for Test 1.}\label{fig3}
\end{figure}

Figure~\ref{fig3} (a) shows the zero level sets of the adaptive
finite element solutions at $t=0.01$, computed by using $\eps=0.01$
and three different tolerances $TOL = 0.01, 0.02$ and $0.04$. The
difference of the three curves is almost invisible, which implies
that we do not need to impose a stringent smallness constraint on
the initial error and the residual (cf. Corollary~\ref{cor2.1}), and
that the continuous dependence estimate of Proposition~\ref{prop2.3}
may be improved.

If we zoom in at the left tip of the curves in Figure~\ref{fig3} (a),
we then find that the distance between the zero level sets for
$TOL=0.04$ and $0.02$ is about $0.00173$, and the distance between
the zero level sets for $TOL=0.02$ and $0.01$ is about $0.0004$ (see
Figure~\ref{fig3} (b)). Since the DOFs at time $0.01$ with respect
to $TOL = 0.01, 0.02$ and $0.04$ are $\mathcal{N}_{0.01} = 12565$,
$\mathcal{N}_{0.02} = 9766$ and $\mathcal{N}_{0.04} = 5995$,
respectively, we have
\[
\frac{1/\mathcal{N}_{0.02}^2-1/\mathcal{N}_{0.01}^2}{1/\mathcal{N}_{0.04}^2
-1/\mathcal{N}_{0.02}^2}\approx 0.2394\approx
0.2312\approx\frac{0.0004}{0.00173}.
\]
Hence, the rate of convergence of the zero level set of the adaptive finite
element solution is about $O(1/\mathcal{N}^2)$. Figure~\ref{fig3}
(c) shows the zero level sets of the adaptive finite element
solution at time $0.01$, computed by using $TOL = 0.02$ and
$\eps=0.08, 0.04, 0.02$ and $0.01$, respectively.


{\bf Test 2:} Consider the Cahn-Hilliard equation
\eqref{e1.1}-\eqref{e1.3} with the initial condition
\begin{equation}\label{E2}
\begin{split}
 u_0(x,y)= &\tanh\bigl(((x-0.3)^2+y^2-0.2^2)/\eps\bigr)
 \tanh\bigl(((x+0.3)^2+y^2-0.2^2)/\eps\bigr)\times\\
 & \tanh\bigl((x^2+(y-0.3)^2-0.2^2)/\eps\bigr)
 \tanh\bigl((x^2+(y+0.3)^2-0.2^2)/\eps\bigr).
\end{split}
\end{equation}

\begin{figure}[hbt]
\centerline{
\includegraphics[scale=0.27]{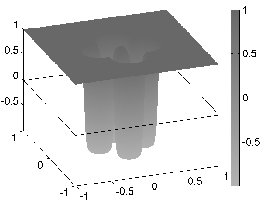}
\includegraphics[scale=0.27]{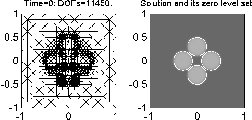}
} \caption{The profile of $u_0$ and its zero level set of Test
2}\label{fig4}
\end{figure}

\begin{figure}[htb]
\centerline{
\includegraphics[scale=0.27]{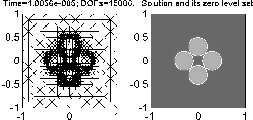}
\includegraphics[scale=0.27]{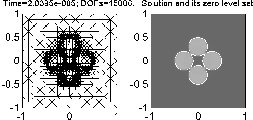}
} \centerline{
\includegraphics[scale=0.27]{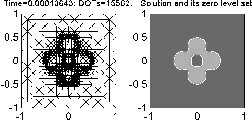}
\includegraphics[scale=0.27]{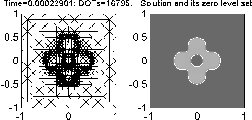}
} \centerline{
\includegraphics[scale=0.27]{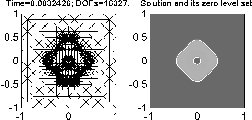}
\includegraphics[scale=0.27]{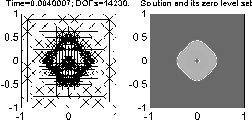}
} \centerline{
\includegraphics[scale=0.27]{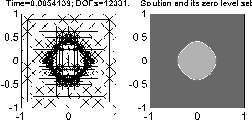}
\includegraphics[scale=0.27]{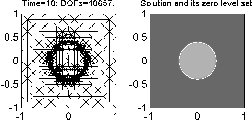}
}  \caption{Snapshots of computed solutions and adaptive meshes for
Test 2}\label{fig5}
\end{figure}

Figure~\ref{fig4} displays the graph of the initial function $u_0$
and its zero level set, which encloses four circles with radius
$0.2$. It also shows the initial mesh and computed initial zero
level set $\Gamma_0^{0.01,h}$. Figure~\ref{fig5} shows snapshots of
the solution (and its zero level set) of the Cahn-Hilliard equation
and the (adaptive) mesh on which the solution is computed at $15$
different time steps. $\eps=0.01$ and $TOL=0.02$ are used in the
simulation. As expected, the fine mesh follows the zero level set as
it moves. We also note that the number of elements in the initial
mesh $\mathcal{T}_0$ is $2520$, the minimum area of the elements is
$1.2207\times 10^{-4}$. If a uniform mesh is used, we need
$\frac{4}{1.2207}\times 10^{4}\approx 32,768$ elements and about
$148,000$ DOFs.


{\bf Test 3:} Consider the Cahn-Hilliard equation
\eqref{e1.1}-\eqref{e1.3} with the following initial condition
\begin{equation}\label{E3}
\begin{split}
 u_0&(x,y)= \tanh\bigl((x^2+y^2-0.15^2)/\eps\bigr)\times\\
 &\tanh\bigl(((x-0.31)^2+y^2-0.15^2)/\eps\bigr)\tanh\bigl(((x+0.31)^2+y^2-0.15^2)/\eps\bigr)\times\\
 & \tanh\bigl((x^2+(y-0.31)^2-0.15^2)/\eps\bigr)\tanh\bigl((x^2+(y+0.31)^2-0.15^2)/\eps\bigr)\times\\
 &\tanh\bigl(((x-0.31)^2+(y-0.31)^2-0.15^2)/\eps\bigr)\times\\
 &\tanh\bigl(((x-0.31)^2+(y+0.31)^2-0.15^2)/\eps\bigr)\times\\
 &\tanh\bigl(((x+0.31)^2+(y-0.31)^2-0.15^2)/\eps\bigr)\times\\
 &\tanh\bigl(((x+0.31)^2+(y+0.31)^2-0.15^2)/\eps\bigr).
\end{split}
\end{equation}

\begin{figure}[hbt]
\centerline{
\includegraphics[scale=0.18]{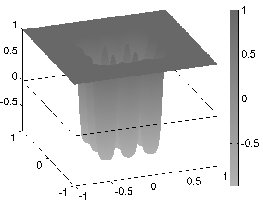}
\includegraphics[scale=0.18]{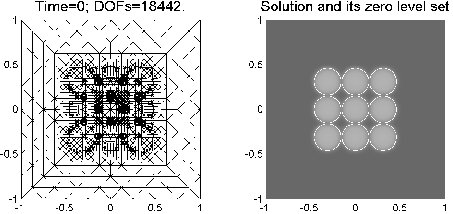}
} \caption{The profile of $u_0$ and its zero level set of Test
3}\label{fig6}
\end{figure}

\begin{figure}[htb]
\centerline{
\includegraphics[scale=0.27]{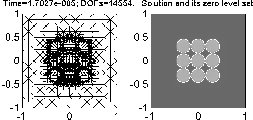}
\includegraphics[scale=0.27]{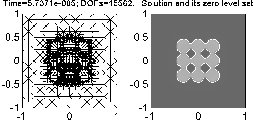}
} \centerline{
\includegraphics[scale=0.27]{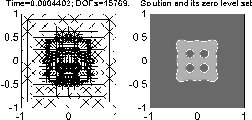}
\includegraphics[scale=0.27]{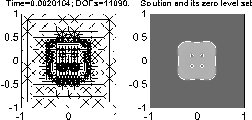}
} \centerline{
\includegraphics[scale=0.27]{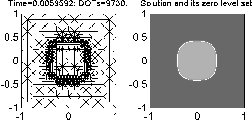}
\includegraphics[scale=0.27]{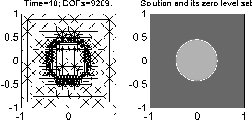}
}  \caption{Snapshots of computed solutions and adaptive meshes for
Test 3}\label{fig7}
\end{figure}

Figure~\ref{fig6} displays the graph of the initial function $u_0$
and its zero level set, which encloses nine circles with radius
$0.15$. It also shows the initial mesh and computed initial zero
level set $\Gamma_0^{0.01,h}$. Figure~\ref{fig7} shows snapshots of
the solution (and its zero level set) of the Cahn-Hilliard equation
and the (adaptive) mesh on which the solution is computed at $15$
different time steps. $\eps=0.01$ and $TOL=0.02$ are used in the
simulation. As expected, the fine mesh follows the zero level set as
it moves. We also note that the number of elements in the initial
mesh $\mathcal{T}_0$ is $4,072$, the minimum area of the elements is
$3.0518\times 10^{-5}$. If a uniform mesh is used, we need
$\frac{4}{3.0518}\times 10^{5}\approx 131,072$ elements and about
$590,000$ DOFs.


%%%%%%%%%%
\bibliographystyle{./abbrv}

\end{document}